\magnification=\magstep1 
\overfullrule=0pt \input amssym.def

\def\eqde{\,{\buildrel \rm def \over =}\,}  
\def\al{\alpha} \def\la{\lambda}  \def\ga{\gamma} 
       \def\i{{\rm i}} 
\def\si{\sigma} \def\eps{\epsilon}    
\def\s{{\cal S}}         
\def\r{\overline{r}}  \def\L{{\Lambda}}  
    \def\D{{\cal D}} \def\A{{\cal A}}
  \def\C{{\Bbb C}}  \def\Z{{\Bbb Z}}
\def\Q{{\Bbb Q}} \def\R{{\Bbb R}}
\font\huge=cmr10 scaled \magstep2
\def\QED{\vrule height6pt width6pt depth0pt}
\font\smcap=cmcsc10

\def\boxit#1{\vbox{\hrule\hbox{\vrule{#1}\vrule}\hrule}}
\def\splus{\,\,{\boxit{$+$}}\,\,}
\def\stimes{\,\,{\boxit{$\times$}}\,\,}
\def\sdot{\,{\boxit{$\cdot$}}\,}

\centerline{{\bf \huge  The Automorphisms of Affine Fusion Rings}}
\bigskip \bigskip   \centerline{{\smcap Terry Gannon}}\medskip
\centerline{{\it Department of 
Mathematical Sciences, University of Alberta,}}
\centerline{{\it  Edmonton, Canada, T6G 2G1}}
\smallskip
\centerline{{e-mail: tgannon@math.ualberta.ca}}
\bigskip\bigskip

\centerline{{\smcap 1. Introduction}}\bigskip 

Verlinde's formula [33]
$$V^{(g)}_{a^1\cdots a^t}=\sum_{b\in \Phi} (S_{0b})^{2(1-g)}{S_{a^1b}\over
S_{0b}}\cdots {S_{a^tb}\over S_{0b}}\eqno(1.1a)$$
arose first in rational conformal field theory (RCFT) as an extremely useful
expression for
the dimensions of conformal blocks on a genus $g$ surface with $t$ punctures.
$\Phi$ here is the finite set of `primary fields'. The matrix $S$ 
comes from a representation of SL$_2(\Z)$ defined by the chiral characters of 
the theory. Contrary to appearances, these numbers $V^{(g)}_{\star\cdots\star}$ will 
always be nonnegative integers.
See the excellent bibliography in [6] for references to the
physics literature.

These numbers are remarkable for also arising in several other contexts:
for example, as  dimensions of spaces of generalised theta
functions; as certain tensor product coefficients in quantum groups and
Hecke algebras at roots of 1 and Chevalley
groups for ${\Bbb F}_p$; as certain knot invariants for 3-manifolds; as
composition laws of the superselection sectors in algebraic quantum field
theories; as dimensions of spaces of intertwiners  in vertex operator
algebras (VOAs); in von Neumann algebras as ``Connes' fusion''; in quantum
cohomology; and in Lusztig's
exotic Fourier transform. See for example  [7,20,19,32,11,10,36,37,26], and
references therein.
 
The more fundamental of these numbers are those corresponding to a sphere
with three punctures. It is more convenient to write these
in the form (called {\it fusion coefficients})
$$N_{ab}^c\eqde V_{a,b,Cc}^{(0)}=\sum_{d\in \Phi}{S_{ad}S_{bd}S^*_{cd}\over S_{0d}}\eqno(
1.1b)$$
where $C$ is a permutation of $\Phi$ called {\it charge-conjugation} and will
be defined below. The fusion
coefficients uniquely determine all other Verlinde dimensions (1.1a). 
 The symmetries of the numbers 
(1.1b), i.e.\ the permutations $\pi$ of $\Phi$ obeying
$$N_{\pi a,\pi b}^{\pi c}=N_{ab}^c\ ,\eqno(1.2)$$
  are precisely the symmetries of all numbers of the form (1.1a).

The point of introducing the $N_{ab}^c$ in (1.1b) is that they define an algebraic
structure, the {\it fusion ring}. Consider all
formal linear combinations of objects $\chi_a$
labelled by the $a\in \Phi$; the multiplication is defined to have structure
constants $N_{ab}^c$:
$$\chi_a\chi_b=\sum_{c\in \Phi}N_{ab}^c\chi_c\eqno(1.3)$$
As an abstract ring, it is not so interesting (the fusion ring over
$\C$ is isomorphic to $\C^{\|\Phi\|}$ with operations defined component-wise;
over $\Q$ it will be a direct sum of number fields).
This is analogous to the character ring of a Lie algebra, which is isomorphic
as a ring to a polynomial ring. Of course it is important in both contexts that 
we have a preferred basis, namely $\{\chi_a\}$, and so proper definitions
of isomorphisms etc.\ must respect
that.

The most important examples of fusion rings are associated to the affine
algebras, and it is to these that this paper is devoted. Their automorphisms
appear explicitly for instance in the classification of modular
invariant (i.e.\ {\it torus})  partition functions
[17,18],
and also in D-branes and boundary conditions for conformal field theory (see
e.g.\ [1]). For instance, fusion-automorphisms (more generally,
-homomorphisms) generate large classes of nonnegative integer
representations of the fusion-ring, each of which is associated to a
boundary ({\it cylinder}) partition function. This will be studied elsewhere.
Also,  whenever the coefficient matrix of the torus partition function
is a permutation matrix (in which case the partition function is called an
{\it automorphism invariant}), we get a fusion ring automorphism. However
most torus partition functions are not automorphism invariants (although
Moore-Seiberg assert that there is a sense in which any torus partition function
can be interpreted as one --- see e.g.\ [3]), and most fusion
ring automorphisms do not correspond to partition functions. Nevertheless,
the two problems are related. The automorphism
invariants for the affine algebras were classified in [17,18]; a Lemma
proved there (our Proposition 4.1 below) involving q-dimensions will
be very useful to us, and  conversely the arguments in Section 4 of this paper 
could be used to considerably simplify the proofs of [17,18].

It is surprising that it is even possible
to find all affine fusion automorphisms, and in fact the arguments turn out to
be rather short. It is remarkable that the answer is so simple: with
few exceptions, they correspond to the Dynkin diagram symmetries.

A related task is determining which affine fusion rings 
 are isomorphic. We answer this in section 5 below; as expected
most fusion rings with different names are nonisomorphic.

\medskip{{\smcap Acknowledgements.}} Most of this paper was written during a
visit to Max-Planck in Bonn, whose hospitality as always was both stimulating
and pleasurable. I also thank Yi-Zhi Huang for clarifying an issue concerning
[8,21].

\vfill\eject\centerline{{\smcap 2. Generalities}} 
\bigskip

\noindent{{\it 2.1. The affine fusion ring}}\medskip

 The source of some of the most interesting 
fusion data are the
affine nontwisted Kac-Moody algebras $X_r^{(1)}$ [23]. Choose any positive 
integer $k$. Consider the (finite) set $P_+=P_+^k(X_r^{(1)})$ of level $k$ integrable highest weights:
$$P_+\eqde \bigl\{\sum_{j=0}^r\la_j\L_j\,
|\,\la_j\in \Z,\ \la_j\ge 0,\ \sum_{j=0}^ra^\vee_j\la_j=k\bigr\}\ ,$$
where $\Lambda_i$ denote the fundamental weights, and  
$a_j^\vee$ are the co-labels, of $X_r^{(1)}$ (the $a_j^\vee$ will be given for
each algebra in \S3).  We will usually drop the (redundant) component
$\la_0\L_0$.
Kac-Peterson [24] found a natural representation of the modular group SL$_2(\Z)$ on
the complex space spanned by the affine characters $\chi_\mu$, $\mu\in P_+$:
most significantly, $\left(\matrix{0&-1\cr1&0}\right)$ is sent to the {\it
Kac-Peterson} matrix $S$ with entries
$$S_{\mu\nu}= c\,
\sum_{w\in \overline{W}}{\rm det}(w)\,\exp[-2\pi \i\,{(w(\mu+\rho)|
\nu+\rho)\over \kappa}]\ .\eqno(2.1a)$$
An explicit expression for the
normalisation constant $c$ is given in e.g.\ [23, Theorem 13.8(a)].
The inner product in (2.1a) is scaled so that the long roots have norm 2.
$\overline{W}$ is the (finite) Weyl group of
$X_r$, and acts on $P_+$ by fixing $\L_0$. The Weyl vector $\rho$ equals
$\sum_i\L_i$, and $\kappa\eqde k+\sum_ia_i^\vee$. This is the matrix
$S$ appearing in (1.1); $\Phi$ there is  $P_+$ here.

The matrix $S$ is symmetric and unitary.              
One of the weights, $k\L_0$, is distinguished and will be denoted `0'.
It is the weight appearing in the denominator of (1.1). A useful fact
is that $$S_{\la 0}>0\qquad {\rm for\ all}\ \la\in P_+\ .$$ 

Equation (2.1a) gives us  the important 
$$\chi_\la[\mu]\eqde {S_{\la\mu}\over S_{0\mu}}=
{\rm ch}_{\overline{\la}}(-2\pi \i
{\overline{\mu}+\overline{\rho}\over \kappa})\ ,\eqno(2.1b)$$
where ${\rm ch}_{\overline{\la}}$ is the Weyl character of the
$X_r$-module ${L}(\overline{\la})$. Together with the
Weyl denominator formula, it provides a useful expression for the
{\it q-dimensions}:
$$\D(\la)\eqde{S_{\la0}\over S_{00}}=\prod_{{\al}>0}
{\sin(\pi \,(\la+\rho\,|\,{\al})/\kappa)\over \sin(\pi\,(\rho\,|\,
{\al})/\kappa)}\ ,\eqno(2.1c)$$
where the product is over the positive roots ${\al}\in
\overline{\Delta}_+$ of $X_r$.
Another consequence of (2.1b) is the  Kac-Walton formula (2.4).

{\it Charge-conjugation} is the order 2 permutation of $P_+$ given by
 $C\la={}^t\la$, the weight contragredient to
$\la$. For instance $C0=0$.
It has the basic property that
$$S_{C\la,\mu}=S_{\la,C\mu}=S_{\la\mu}^*\eqno(2.2a)$$
 and $S^2=C$.
$C$ corresponds to a symmetry of the (unextended) Dynkin diagram of $X_r$,
as we will see next section. 

Related to $C$ are all the other symmetries of the unextended Dynkin diagram.
 We call these {\it conjugations}. The only $X_{r}^{(1)}$ with nontrivial
 conjugations other than charge-conjugation are $D_{even}^{(1)}$.

Another important symmetry of the matrix $S$ is called {\it simple-currents}.
Any weight $j\in P_+$ with {q-dimension} $\D(j)=1$,
is called a simple-current. To any such weight $j$ is associated a permutation
$J$  of $P_+$ and a function $Q_j:P_+\rightarrow\Q$, such that $J0=j$ and
$$S_{J\la,\mu}=\exp[2\pi\i\, Q_j(\mu)]\,S_{\la\mu}\eqno(2.2b)$$
The simple-currents form an abelian group, given by composition of the
permutations $J$.

All simple-currents for the affine algebras were classified in [12], and with
one unimportant exception ($E_8^{(1)}$ at level 2) correspond to symmetries
of the extended Coxeter--Dynkin diagram of $X_r^{(1)}$. The simplest
proof would use the methods of Proposition 4.1 below. For a more intrinsically
algebraic interpretation of these simple-currents, see [25] where their
group is denoted $W_0^+$.

Evaluating $S_{J\la,j'}$ in two ways gives the useful
$$Q_{j'}(J\la)\equiv Q_j(j')+Q_{j'}(\la)\qquad({\rm mod}\ 1)\eqno(2.2c)$$
and hence the reciprocity $Q_j(j')=Q_{j'}(j)$.

For each $X_r$, the inner products $(\la|\mu)$ of weights are rational;
let $N$ denote the least common denominator. E.g.\ for $A_r$ this is
$N=r+1$, while for $E_8$ it is $N=1$. Choose any integer $\ell$ coprime
to $\kappa N$. Then for any $\la\in P_+$ there is a unique weight $\la^{(\ell)}
\in P_+$, coroot $\alpha$, and (finite) Weyl element $\omega$ such that
$$\ell\,(\la+\rho)=\omega(\la^{(\ell)}+\rho)+\kappa \alpha\ .$$
This is simply the statement that the affine Weyl orbit of $\ell\,(\la+\rho)$
intersects the set $P_++\rho$ at precisely one point (namely $\la^{(\ell)}+\rho$).
Write $\eps'_\ell(\la)={\rm det}\,\omega=\pm 1$. Then [16]
$$\eps'_\ell(\la)\,S_{\la^{(\ell)},\mu}=\eps'_\ell(\mu)\,S_{\la,\mu^{(\ell)}}
\eqno(2.3a)$$
This has an obvious interpretation as a Galois automorphism [4]: the field
generated over $\Q$  by all entries $S_{\la\mu}$ lies in the
cyclotomic field $\Q[\xi_{4N\kappa}]$ where $\xi_n$ denotes the root of unity
$\exp[2\pi\i/n]$; for any $\si_\ell\in{\rm Gal}
(\Q[\xi_{4N\kappa}]/\Q)\cong\Z_{4N\kappa}^\times$, there will be
a function $\eps_\ell:P_+\rightarrow\{\pm 1\}$ such that
$$\si_\ell(S_{\la\mu})=\eps_\ell(\la)\,S_{\la^{(\ell)},\mu}=\eps_\ell(\mu)\,
S_{\la,\mu^{(\ell)}}\ .\eqno(2.3b)$$
$\eps_\ell(\la)/\eps'_\ell(\la)=\si_\ell(c)/c$ is an unimportant
sign independent of $\la$. This Galois action will play a fairly important
role in this paper. Note that $\si_{-1}=C$, so this action can be thought of as a
generalisation of charge-conjugation. Note also that 
$\si_\ell\circ J=J^\ell\circ\si_\ell$. 

The fusion coefficients (1.1b)
are usually computed by the Kac-Walton formula [23 p.\ 288, 35] (there are other
codiscoverers) in terms of the tensor product 
multiplicities $T_{{\la}{\mu}}^{{\nu}}\eqde\,$mult$_{
L(\overline{\la})\otimes L(\overline{\mu})}(L(\overline{\nu}))$ 
in $X_r$:
$$N_{\la\mu}^\nu=\sum_{w\in{W}}{\rm det}(w)\,T_{{\la}{\mu}}^{
{w.\nu}}\ ,\eqno(2.4)$$
where $w.\ga\eqde w(\ga+\rho)- \rho$ and $W$ is the affine Weyl group of
$X_r^{(1)}$ (the dependence of $N_{\la\mu}^\nu$ on $k$ arises through the
action of $W$). We shall see shortly that these fusion coefficients, now
manifestly integral,  are in fact nonnegative.
Let ${\cal R}(X_{r,k})$ denote the corresponding fusion ring.

A handy consequence of
(2.4) that whenever $k$ is large enough that $\la+\mu\in
P_+^k(X_r^{(1)})$ (i.e.\ that $\sum_{i=1}^ra_i^\vee(\la_i+\mu_i)\le k$),
then $N_{\la\mu}^\nu=T_{\la\mu}^\nu$.

It will sometimes be convenient to collect these coefficients in matrix
form as the {\it fusion matrices} $N_\la$, defined by $(N_\la)_{\mu\nu}
=N_{\la\mu}^\nu$. For instance, $N_0=I$ and,
more generally, $N_j$ is the permutation matrix associated to $J$.

The importance of (charge-)conjugation and simple-currents for us is that
they respect fusions:
$$\eqalignno{N_{C\la,C\mu}^{C\nu}=&\,N_{\la\mu}^\nu&(2.5a)\cr
N_{J\la,J'\mu}^{JJ'\nu}=&\,N_{\la\mu}^{\nu}&(2.5b)\cr
N_{\la\mu}^\nu\ne 0\ \Rightarrow\ Q_j(\la)+&\,Q_j(\mu)\equiv Q_j(\nu)\qquad
{(\rm mod}\ 1)&(2.5c)\cr}$$
for any simple-currents $J,J',j$.

For example, for  ${\cal R}(A_{1,k})$ we may take $P_+=\{0,1,\ldots,k\}$
(the value of $\la_1$),
and then the Kac-Peterson matrix is $S_{ab}=\sqrt{{2\over k+2}}\,\sin(\pi\,{(a+1)\,
(b+1)\over k+2})$. Charge-conjugation $C$ is trivial here, but $j=k$ is a
simple-current corresponding to permutation $Ja=k-a$ and function
$Q_j(a)=a/2$. The Galois action sends $a$ to the unique weight $a^{(\ell)}\in
P_+$ satisfying $a^{(\ell)}+1\equiv \pm\ell\,(a+1)$ (mod $2k+4$), where
that sign there equals $\i^{\ell-1}\eps'_\ell(a)$. The fusion coefficients are given by
$$N_{ab}^c=\left\{\matrix{1&{\rm if}\ c\equiv a\!+\!b
\ ({\rm mod\ 2)\ and}\ |a\!-\!b|\le c\le{\rm min}\{a\!+\!b,
2k\!-\!a\!-\!b\}\cr0&{\rm otherwise}\cr}\right.\ .$$

Equation (2.4) tells us the affine fusion rules are the structure constants for the
ring Ch$(X_r)/{\cal J}_k$ where Ch$(X_r)$ is the character ring for all
finite-dimensional $X_r$-modules, and ${\cal J}_k$ is the subspace spanned
by the elements ch$_{\overline{\mu}}-({\rm det}\,w){\rm ch}_{\overline{w.\mu}}$.
Finkelberg [8] proved that this ring is isomorphic to the K-ring of
a ``sub-quotient'' $\widetilde{\cal O}_{k}$ of Kazhdan-Lusztig's category of
level $k$
integrable highest weight $X_r^{(1)}$-modules, and to Gelfand-Kazhdan's
category $\widetilde{\cal O}_q$ coming from finite-dimensional modules of the
quantum group $U_qX_r$ specialised
to the root of unity $q=\xi_{2m\kappa}$ for appropriate choice of $m\in
\{1,2,3\}$. They also arise from the Huang-Lepowsky coproduct [21] for the modules
of the VOA $L(k,0)$.  Because of these isomorphisms, we get that
the $N_{\la\mu}^\nu$ do indeed lie in $\Z_\ge$, for any affine algebra.

A useful way of identifying weights in affine Weyl orbits involves computing
q-dimensions and norms. Q-dimensions vary by at most a sign
while norms are constant mod $2\kappa$:
$\D(w.\la)={\rm det}\,(w)\,\D(\la)$ and $(w\la|w\la)\equiv(\la|\la)$ (mod $2\kappa$).
The point is that for exceptional algebras at small levels, the highest weights
can often be distinguished by the pair $(\D(\la),(\la+\rho|\la+\rho)\ ({\rm
mod}\ 2\kappa))$. For example this is true of $E_{8,5},E_{8,6},F_{4,4}$.
This is a useful way in practise to use both (2.4) and the Galois action (2.3).

An  important
property obeyed by the matrix $S$ for any classical algebra $X_r$ is 
{\it rank-level duality}. The first appearance of this curious duality 
seems to be by Frenkel [9], but by now many aspects and generalisations have
been explored in the literature.
For $A_r^{(1)}$, it is related to the existence of mutually
commutative affine subalgbras $\widehat{{\rm sl}(n)}$ and $\widehat{{\rm sl}
(k)}$ in $\widehat{{\rm gl}(nk)}$.  Witten has another interpretation of
it [37]: he found a natural map (a ring homomorphism) from the quantum
cohomology of the Grassmannian $G(k,N)$, to the fusion ring of the
algebra ${\rm u}(k)\cong {\rm su}(k)\oplus {\rm u}(1)$ at level $(N-k,N)$.
Witten used the
duality between $G(k,N)$ and $G(N-k,N)$ to show that the fusion rings of
u$(k)$ level $(N-k,N)$ and u$(N-k)$ level $(k,N)$ should coincide.
 A considerable generalisation, applying to any VOA (or RCFT),
has been conjectured by Nahm [30], and relates to the  natural involution
$\sum_i[x_i]\leftrightarrow\sum_i[1-x_i]$ of torsion elements of the Bloch
group.

The Kac-Peterson matrices of $\widehat{{\rm sl}(\ell)}$ level $k$
and $\widehat{{\rm sl}(k)}$ level $\ell$ are related, as are those of $C_{r,k}$
and $C_{k,r}$, and $\widehat{{\rm so}(\ell)}$ level $k$ and $\widehat{{\rm so}(k)}$
 level $\ell$. We will need only the symplectic one; the
details will be given in \S3.3.

\bigskip\noindent{{\it 2.2. Symmetries of fusion coefficients}}\medskip

{{\smcap Definition 2.1.}}\quad {\it
By an {\it isomorphism}  between fusion rings ${\cal R}(X_{r,k})$
and ${\cal R}(Y_{s,m})$ (with fusion coefficients $N$ and $M$ respectively)
we mean a bijection $\pi:P_+^k(X_{r}^{(1)})\rightarrow
P_+^m(Y_{s}^{(1)})$ such that
$$N_{\la,\mu}^\nu=M_{\pi\la,\pi\mu}^{\pi\nu}\qquad\forall\,\la,\mu,\nu\in
P_+(X_{r,k})\ .\eqno(2.6)$$
When $X_{r,k}=Y_{s,m}$ we call $\pi$ an automorphism or}
fusion-symmetry. {\it Call the pair
of permutations $\pi,\pi'$ an} $S$-symmetry {\it if}
$$S_{\pi \la,\pi' \mu}=S_{\la\mu}\qquad \forall \la,\mu\in P_+\ .$$

The lemma below tells us that fusion- and $S$-symmetries form two isomorphic groups; the former
 we will label ${\cal A}(X_{r,k})$. Equation
(2.5a) says that the charge-conjugation $C$, and more generally any conjugation,
 is a fusion-symmetry, while (2.2a) says
$(C,C)$ is an $S$-symmetry. Because $N_0=I=M_{\tilde{0}}$, $N_{\la\mu}^0=C_{\la\mu}$
and $M_{\tilde{\la},\tilde{\mu}}^{\tilde{0}}=\widetilde{C}_{\tilde{\la},
\tilde{\mu}}$ (we use tilde's to denote quantities in $Y_{s}^{(1)}$ level $m$), any
isomorphism $\pi$ must obey $\pi 0=\tilde{0}$ and $\widetilde{C}\circ\pi=
\pi\circ C$.
More generally, since $N_\la$ is a permutation matrix of order $n$ iff
$\la$ is a simple-current of order $n$, we see that an isomorphism
sends simple-currents to simple-currents of equal order. We get
 $$\pi(J\mu)=\pi(j)\,\pi(\mu)\ .\eqno(2.7a)$$
For instance $\pi$ must send $J$-fixed-points to $\pi(J)$-fixed-points.

More generally, a {\it fusion-homomorphism} $\pi$
is defined in the obvious algebraic way. It turns out that for such a
$\pi$, $\pi\la=\pi\mu$ iff $\mu=J\la$ for some simple-current $J$
for which $\pi(J0)=\tilde{0}$. Moreover, $\pi(J0)=\tilde{0}$ is possible
only if there are no $J$-fixed-points. When $\pi$ is one-to-one (e.g.\
when there are no nontrivial simple-currents in $P_+^k(X_r^{(1)})$),
then $\pi$ obeys (2.6). Fusion-homomorphisms will be studied elsewhere.

The key to finding
fusion-symmetries is the following Lemma.

\medskip{{\smcap Lemma 2.2.}}\quad {\it Let $\widetilde{S}$ be the Kac-Peterson
matrix for $Y_{s}^{(1)}$ level $m$. Then a bijection $\pi:P_+^k(X_r^{(1)})
\rightarrow P_+^m(Y_s^{(1)})$ defines an isomorphism of fusion rings iff
there exists some bijection $\pi':P_+^k(X_r^{(1)})\rightarrow
P_+^m(Y_s^{(1)})$ such that $S_{\la\mu}=\widetilde{S}_{\pi\la,\pi'\mu}$
for all $\la,\mu\in P_+^k(X_r^{(1)})$. In particular, a permutation $\pi$ is a
fusion-symmetry iff $(\pi,\pi')$ is an $S$-symmetry for some $\pi'$.}\medskip 

\noindent{{\it Proof.}}\quad The equality $N_{\la\mu}^\nu=M_{\pi \la,\pi \mu}^{\pi \nu}$
means that, for each $\mu$, the column vectors $({\underline x}_\mu)_\nu=
\widetilde{S}_{\pi \nu,\pi \mu}$ are
simultaneous eigenvectors for the fusion matrices $N_\la$, with eigenvalues
$\widetilde{S}_{\pi \la,\pi \mu}/
\widetilde{S}_{0,\pi \mu}$. It is easy to see from Verlinde's formula (1.1b) that any simultaneous
eigenvector for all fusion matrices must be a scalar multiple of some
column of $S$. Thus there must be  a permutation $\pi''$ of $P_+^k(X_r^{(1)})$
and scalars
$\alpha(\mu)$ such that $\widetilde{S}_{\pi \nu,\pi \mu}=\alpha(\mu)\,S_{\nu,\pi''\mu}$. Taking $\nu=0$
forces $\alpha(\mu)>0$, and then unitarity forces $\alpha(\mu)=1$. \qquad \QED
\medskip

Let $\pi$ be any isomorphism, and let $\pi'$ be as in the Lemma. Then $\pi'$
is also an isomorphism, with $(\pi')'=\pi$. Equation
(2.2b) implies for all $\la\in P_+$ and all simple-currents $j$, that
$$Q_j(\la)\equiv \widetilde{Q}_{\pi'j}(\pi\la)\equiv \widetilde{Q}_{\pi j}(
\pi'\la)\qquad ({\rm mod}\ 1)\ .\eqno(2.7b)$$

Another quick consequence of the Lemma is that for any Galois automorphism $\si_\ell$
and isomorphism $\pi$,
we have $\tilde{\eps}_\ell(\pi\la)=\eps_\ell(\la)$ and $\pi(\la^{(\ell)})=(\pi\la)^{(\ell)}$.
To see this, apply the invertibility of $S$ to the equation
$$\eps_\ell(\la)\,S_{\la^{(\ell)},\mu}=\si_\ell S_{\la\mu} =
\si_\ell \widetilde{S}_{\pi\la,\pi'\mu}=\tilde{\eps}_\ell(\pi\la)\,
\widetilde{S}_{(\pi\la)^{(\ell)},\pi'\mu}=\tilde{\eps}_\ell(\pi\la)\,
{S}_{\pi^{-1}(\pi\la)^{(\ell)},\mu}\ .$$

A very useful notion for studying the fusion ring  is that of {\it
fusion-generator}, i.e.\ a subset $\Gamma=\{\ga_1,\ldots,\ga_m\}$ of $P_+$
which generates ${\cal R}(X_{r,k})$ as a ring.  Diagonalising, this is equivalent
to requiring that there are $m$-variable polynomials $P_\la(x_1,\ldots,x_m)$
such that
$${S_{\la\mu}\over S_{0\mu}}=P_\la({S_{\ga_1\mu}\over S_{0\mu}},\ldots,{S_{\ga_m
\mu}\over S_{0\mu}})\qquad \forall \la,\mu\in P_+\ .$$
Let $(\pi,\pi')$ be an $S$-symmetry, and suppose we know that $\pi \ga=\ga$ for
all $\ga$ in the fusion-generator $\Gamma$. Then for any $\la\in P_+$,
$${S_{\la\mu}\over S_{0\mu}}={S_{\pi \la,\pi' \mu}\over S_{0,\pi' \mu}}=P_{\pi
\la}({S_{\ga_1,\pi'\mu}\over S_{0,\pi'\mu}},\ldots)=P_{\pi \la}({S_{\ga_1\mu}
\over S_{0\mu}},\ldots)={S_{\pi \la,\mu}\over S_{0\mu}}$$
for all $\mu\in P_+$, so $\pi \la=\la$. 

One of the reasons fusion-symmetries for the affine algebras are so tractible
is the existence of small fusion-generators. In particular, because we know
that any Lie character ch$_{\overline{\mu}}$ for $X_r$ can be written as a polynomial in the
fundamental characters ${\rm ch}_{\L_1},\ldots,{\rm ch}_{\L_r}$, we know from
(2.1b) that $\Gamma=
\{\L_1,\ldots,\L_r\}$ is a fusion-generator for $X_r^{(1)}$ at any level $k$
sufficiently large that $P_+$ contains all $\L_i$ (in other
words, for any $k\ge {\rm max}_i\,a_i^\vee$). In fact, it is easy to show
[18] that a fusion-generator valid for any $X_{r,k}$ is $\{\L_1,\ldots,\L_r\}
\cap P_+$. Smaller fusion-generators usually
exist --- for example $\{\L_1\}$ is
a fusion-generator for $A_{8,k}$ whenever $k$ is even and coprime to 3.

\bigskip
\noindent{{\it 2.3. Standard constructions of fusion-symmetries}}
\medskip

Simple-currents are a large source of fusion-symmetries.
Let $j$ be any simple-current of order $n$. Choose any number $a\in\{0,1,
\ldots,n-1\}$ such that $${\rm gcd}(naQ_j(j)+1,n)=1\ .$$ Any solution to this defines
a fusion-symmetry $\la\mapsto J^{naQ_j(\la)}\la$, which we shall denote
$\pi[a]$ or $\pi_j[a]$. Note that from (2.2b), (2.5b) and (2.5c) that any
$\pi=\pi[a]$, $a\in\Z$, obeys the relation $N_{\pi\la,\pi\mu}^{\pi\nu}=\
N_{\la\mu}^\nu$ when $N_{\la\mu}^\nu\ne 0$ (it  would in fact be a
fusion-endomorphism --- see \S2.2); the `gcd' condition forces $\pi[a]$ to be a {\it permutation}.
Choosing $b\equiv-a\,(naQ_j(j)+1)^{-1}$ (mod $n$), we find that $(\pi[a],\pi[b])$
is an $S$-symmetry.

When the group of simple-currents is not cyclic,
this construction can be generalised in a natural way, and the resulting
fusion-symmetry will be parametrised by a matrix $(a_{ij})$.
We will meet these in \S3.4.

We will call these {\it simple-current automorphisms}. The first
examples of these were found by Bernard [2], and were generalised further in [31].

For any affine algebra $X_r^{(1)}$ and any sufficiently high level, we will see
in the next section that its fusion-symmetries consist entirely of
simple-current automorphisms and conjugations. For this reason, any other
fusion-symmetry is called {\it exceptional}.

There is another general construction of fusion-symmetries, generalising $C$,
although it
yields few new examples for the affine fusion rings.
If the Galois automorphism $\si_\ell$ is such that $0^{(\ell)}$ is a
simple-current $j$ --- equivalently, that $\si_\ell(S_{00}^2)=S_{00}^2$ ---
 then the permutation
$$\pi\{\ell\}:\la\mapsto J(\la^{(\ell)})$$
is a fusion-symmetry. The simplest example is $\pi\{-1\}=C$.
We call $\pi\{\ell\}$ a {\it Galois fusion-symmetry}. A special case of these
 was  given in [13]. To see that $\pi\{\ell\}$ works, note from
$$\eps_\ell(\la)\,S_{\la^{(\ell)},0}=\si_\ell S_{\la 0} =\eps_\ell(0)\,
e^{2\pi\i\,Q_j(\la)}S_{\la 0}$$
that $\eps_\ell(\la)\,\eps_\ell(0)=e^{2\pi\i\,Q_j(\la)}$. Hence
$$S_{J\la^{(\ell)},\mu}=e^{2\pi\i\,Q_j(\mu)}\eps_\ell(\la)\,\si_\ell(S_{\la\mu})
=e^{2\pi\i\,Q_j(\mu)}\,\eps_\ell(\la)\,\eps_\ell(\mu)\,S_{\la,\mu^{(\ell)}}
=S_{\la,J\mu^{(\ell)}}$$
and so $(\pi\{\ell\},\pi\{\ell\}^{-1})$ is an $S$-symmetry. Incidentally,
$J$ will always be order 1 or 2 because $2\,Q_j(\la)\in\Z$ for all
$\la\in P_+$.

Simple-currents (2.2), the Galois action (2.3), and the corresponding
fusion-symmetries have analogues in arbitrary (i.e.\ not necessarily
affine) fusion rings.

\bigskip\bigskip\centerline{{\smcap 3. Data for the Affine Algebras.}}\bigskip

Our main task in this paper is to find and construct
all fusion-symmetries for the affine algebras
$X_r^{(1)}$, for simple $X_r$. In this section we state the results,
 and in the next section we prove the completeness of our
lists. Recall the simple-current automorphism $\pi[a]$ and Galois
automorphism $\pi\{\ell\}$ defined in \S2.3, and the notation $\kappa=k+h^\vee$.
It will be convenient to write `$X_{r,k}$' for `$X_r^{(1)}$ and level
$k$'. We write ${\cal S}$ for the group
of symmetries of the extended Dynkin diagram.

\bigskip\noindent{\it 3.1. The algebra} $A_r^{(1)}$, $r\ge 1$

\medskip Define $\r=r+1$ and $n=k+\r$.
The level $k$ highest weights of $A_r^{(1)}$ constitute the set 
$P_{+}$ of
$\r$-tuples $\la=(\la_0,\ldots,\la_r)$ of non-negative integers 
obeying $\sum_{i=0}^r\la_i=k$. The  Dynkin diagram symmetries form the
dihedral group $\s={\frak D}_{r+1}$; it
is generated by the charge-conjugation $C$ 
and  simple-current $J$ given by  $C\la=(\la_0,\la_r,\la_{r-1},\ldots,\la_1)$
and $J\la=(\la_r,\la_0,\la_1,\ldots,\la_{r-1})$, with $Q_{J^a}(\la)=a\,t(\la)/\r$
for $t(\la)\eqde \sum_{j=1}^r j\la_j$. Note that $C=id.$ for $A_1^{(1)}$. 

The Kac-Peterson relation (2.1b) for $A_{r,k}$ takes the form
$${S_{\la\mu}\over S_{0\mu}}=\exp[-2\pi\i{t(\la)\,t(\mu)\over \kappa\,\r}]\,\,
s_{(\la)}
(\exp[-2\pi\i{(\mu+\rho)(1)\over \kappa},\ldots,\exp[-2\pi\i{(\mu+\rho)(\r)
\over\kappa}])\ ,\eqno(3.1)$$
where $s_{(\la)}(x_1,\ldots,x_{r+1})$ is the Schur polynomial (see e.g.\
[27]) corresponding to the partition $(\la(1),\ldots,\la(\r))$, and
where $\nu(\ell)
=\sum_{i=\ell}^r\nu_i$ for any weight $\nu$. In other words,
$S_{\la\mu}/S_{0\mu}$ is the Schur polynomial corresponding to $\la$,
evaluated at roots of 1 determined by $\mu$.

The fusion (derived from the Pieri rule and (2.4))
$$\L_1\stimes \L_\ell=\L_{\ell+1}\splus(\L_1+\L_\ell)\ ,$$
valid for $k\ge 2$ and $1\le \ell<r$, will be useful.

There are no exceptional fusion-symmetries for $A_r^{(1)}$:

\medskip{\smcap Theorem 3.A.} {\it The fusion-symmetries for $A_r^{(1)}$ level
$k$ are $C^i\pi[a]$, for $i\in\{0,1\}$ and any integer $0\le a\le r$
for which $1+ka$ is coprime to $r+1$.}\medskip

To avoid redundancies in the Theorem, for $r=1$ or $k=1$ take $i=0$ only.
If we write $\r=r'r''$, where $r'$ is coprime to $k$ and $r''|k^\infty$,
then the number of simple-current automorphisms will exactly equal
$r''\cdot\varphi(r')$, where $\varphi$ is the Euler totient. The $\pi[a]$
commute with each other, and with $C$.

For example, for $A_{1,k}$ when $k$ is odd, there is no nontrivial
fusion-symmetry. When $k$ is even, there is exactly one, sending
$\la=\la_1\L_1$ to $\la$ (for $\la_1$ even) or $J\la=(k-\la_1)\L_1$
(for $\la_1$ odd). For $A_{2,k}$, there are either six or four 
fusion-symmetries, depending on whether or not 3 divides $k$.

\vfill\eject\noindent{\it 3.2. The algebra} $B_r^{(1)}$, $r\ge 3$

\medskip
A weight $\la$ in $P_+$ satisfies $k=\la_0+\la_1+2\la_2+\cdots+2\la_{r-1}
+\la_r$, and $\kappa=k+2r-1$.
The charge-conjugation is trivial, but there is a
simple-current: $J\la=(\la_1,\la_0,
\la_2,\ldots,\la_r)$. It has $Q(\la)=\la_r/2$.  

The only fusion products we need are
$$\eqalign{\L_1\stimes\L_i=&\,\L_{i-1}\splus\L_{i+1}\splus(\L_1+\L_i)\cr
\L_1\stimes(\ell\L_r)=&\,(\ell\L_r)\splus(\L_1+\ell\L_r)\splus(\L_{r-1}+
(\ell-2)\L_r)\cr} $$
for all $1\le i<r-1$, $k>2$, and $0<\ell<k$, where we drop
`$\L_{r-1}+(\ell-2)\L_r$' if $\ell=1$. We will also use the character formula (2.1b)
$$\chi_{\L_1}[\la]={S_{\L_1\la}\over S_{0\la}}=2\sum_{\ell=1}^r\cos(2\pi{\la^+
(\ell)\over \kappa})+1 \ ,\eqno(3.2)$$
where $\la^+(\ell)=(\la+\rho)(\ell)$ and 
$$\la(\ell)=\sum_{i=\ell}^{r-1}\la_i+{1\over 2}\la_r\ .$$       

For $k=2$ ($\kappa=2r+1$) there are several Galois fusion-symmetries --- one for each
Galois automorphism, since $S_{00}^2={1\over 4\kappa}$ is rational. In 
particular, define $\ga^i=\ga^{\kappa-i}=\L_i$ for $i=1,2,\ldots,r-1$,
and $\ga^r=\ga^{r+1}=2\L_r$. Then for any $m$ coprime to $\kappa$, $\pi\{m\}$ fixes
$0$ and $J$, sends $\ga^a$ to $\ga^{ma}$ (where the superscript is taken mod 
$\kappa$), and stabilises $\{\L_r,J\L_r\}$ ($\pi\{m\}\L_r=\L_r$ iff the Jacobi
symbol $({\kappa\over m})$ equals +1).

Why is $k=2$ so special here? One reason is that rank-level duality associates
$B_{r,2}$ with u$(1)_{2r+1}$, and it is easy to confirm that $\widehat{{\rm
u}(1)}$
has a rich variety of fusion-symmetries (and modular invariants) coming from its
simple-currents. Also, the $B_{r,2}$ matrix $S$ formally looks like the
character table of the dihedral group and for some $r$ actually equals
the Kac-Peterson matrix $S$ associated to the dihedral group ${\frak D}_{\sqrt{\kappa}}$
twisted by an appropriate 3-cocycle [5] --- finite group modular data
tends to have significantly more modular invariants and fusion-symmetries
than e.g.\ affine modular data.\medskip

{\smcap Theorem 3.B.} {\it The fusion-symmetries of $B_r^{(1)}$ level $k$
for $k\ne 2$ are
$\pi[1]^i$ where $i\in\{0,1\}$. For $k=2$ a fusion-symmetry will equal $\pi[1]^i
\,\pi\{m\}$ for $i\in\{0,1\}$ and $m\in\Z_\kappa^{\times}$, $1\le m\le r$.}

\medskip When $k=1$, $\pi[1]$ is trivial. We have ${\cal
F}(B_{r,2})\cong \Z_2\times(\Z_{2r+1}^\times/\{\pm 1\})$.

\bigskip\noindent{\it 3.3. The algebra} $C_r^{(1)}$, $r\ge 2$

\medskip
A weight $\la$ of $P_+$ satisfies $k=\la_0+\la_1+\cdots+\la_r$ and $\kappa=k+r+1$. 
Charge-conjugation $C$ again is trivial, and there is a simple-current $J$
defined by $J\la=(\la_r,\la_{r-1},\ldots,\la_1,\la_0)$, with
$Q(\la)=(\sum_{j=1}^r \, j\la_j)/2$. 

Choose any $\la\in P_+$. The Young
diagram for $\la$ is defined in the usual way:  for $1\le \ell\le r$, the
$\ell$th row  consists of $\la(\ell)\eqde\sum_{i=\ell}^r\la_i$ boxes.
Let $\tau\la$ denote the $C_{k,r}$ weight whose diagram is
the transpose of that for $\la$. (For this purpose the algebra $C_1$ may be
identified with $A_1$.) For example, $\tau\L_a=
a\tilde{\L}_1$, where we use tilde's to denote the quantities of 
$C_{k,r}$. In fact, $\tau:P_+(C_{r,k})\rightarrow P_+(C_{k,r})$ is a bijection.
Then
$$\tilde{S}_{\tau\la,\tau\mu}=S_{\la\mu}\ .$$

This rank-level duality for $C_r^{(1)}$ is especially interesting, as it defines
a fusion ring isomorphism ${\cal R}(C_{r,k})\cong{\cal R}(C_{k,r})$ (see
\S5). When $k=r$, we get a nontrivial fusion-symmetry: 
$\pi_{{\rm rld}}\la\eqde\tau\la$.

The only fusion product we need is
$$\L_1\stimes\L_i=\L_{i-1}\splus\L_{i+1}\splus(\L_1+\L_i)\ ,$$
valid for $i< r$ and $k\ge 2$. The following character formula (2.1b) will also
be used:
$$\chi_{\L_1}[\la]={S_{\L_1\la}\over S_{0\la}}=2\sum_{\ell=1}^r\cos(\pi{\la^+
(\ell)\over \kappa})\ ,\eqno(3.3)$$
where $\la^+(\ell)=(\la+\rho)(\ell)$ as before.

\medskip{\smcap Theorem 3.C.} {\it The fusion-symmetries for $C_r^{(1)}$ level 
$k$, when $k\ne r$ and either $k$ or $r$ is even, are $\pi[1]^i$ for $i\in\{0,1\}$.
When $k\ne r$ but both $k$ and $r$ are odd, then there is no nontrivial
fusion-symmetry. When $k=r$, they are
$\pi[1]^i\,\pi_{{\rm rld}}^j$ ($k$ even) or $\pi[1]^i$ ($k$ odd),
for $i,j\in\{0,1\}$.}

\medskip When $r=k$ is even, ${\cal A}(C_{r,k})\cong \Z_2\times\Z_2$.

\bigskip\noindent{\it 3.4. The algebra} $D_{r}^{(1)}$, $r\ge 4$

\medskip
A weight $\la$ of $P_+$ satisfies $k=\la_0+\la_1+2\la_2+\cdots+2\la_{r-2}+
\la_{r-1}+\la_r$, and $\kappa=k+2r-2$. 
For any $r$, there are the conjugations $C_0=id.$ and $C_1\la = (\la_0,\la_1,
\ldots, \la_{r-2},\la_r,\la_{r-1})$.
The charge-conjugation $C$ equals $C_1$ for odd $r$, and $C_0$ for even $r$.
When $r=4$ there are four additional conjugations; these six $C_i$ correspond
to all permutations of the $D_4^{(1)}$ Dynkin labels $\la_1,\la_3,\la_4$.

There are three non-trivial simple-currents, $J_v$, $J_s$ and $J_v
J_s$. Explicitly, we have
$J_v\la=(\la_1,\la_0,\la_2,\ldots,\la_{r-2},\la_r,\la_{r-1})$ with
$Q_v(\la)=(\la_{r-1}+\la_r)/2$, and
$$J_s\la=\cases{
(\la_r,\la_{r-1},\la_{r-2},\ldots,\la_1,\la_0) & if $r$ is even, \cr
(\la_{r-1},\la_r,\la_{r-2},\ldots,\la_1,\la_0) & if $r$ is odd,
\cr}$$
with $Q_s(\la)=(2\sum_{j=1}^{r-2}j\la_j-(r-2)
\la_{r-1}-r\la_r)/4$.  From this we compute $Q_s(J_s0)=-rk/4$.

The fusion products we need  are 
$$\eqalignno{\L_1\stimes\L_i=&\,\L_{i-1}\splus \L_{i+1}\splus(\L_1+\L_i)\ &\cr
\L_1\stimes \L_r=&\,\L_{r-1}\splus (\L_1+\L_r)\ ,&\cr}$$
valid for all $1\le i< r-2$ and $k>2$. We also will use the character
formula  (2.1b)
$$\chi_{\L_1}[\la]={S_{\L_1\la}\over S_{0\la}}=2\sum_{\ell=1}^r\cos(2\pi{
\la^+(\ell)\over \kappa})\ ,\eqno(3.4)$$
where $\la^+(\ell)=(\la+\rho)(\ell)$ and the orthonormal components
$\la(\ell)$ are defined by
$\la(\ell)=\sum_{i=\ell}^{r-1}\la_i+{\la_{r}-\la_{r-1}\over 2}$.

The simple-current automorphisms are as follows, and depend on whether $r$ and $k$
are even or odd. When $r$ is odd, the group of simple-currents
is generated by $J_s$. If in addition $k$ is odd,
there will be only two simple-current automorphisms: $\pi=\pi'=\pi[a]=J_s^{4aQ_s}$ for $a\in\{0,2\}$.
If instead $k$ is even, there will be four simple-current automorphisms: $\pi=\pi[a]$
and $\pi'=\pi[ak-a]$ for $0\le a\le 3$. When $k\equiv 2$ (mod 4),
these form the group $\Z_2\times \Z_2$, otherwise when $4|k$ the group
is $\Z_4$.

When $r$ is even, the simple-currents are  generated by both $J_v$ and $J_s$.
If in addition $k$ is even, we have 16 simple-current automorphisms: 
$$\pi=\pi\left[\matrix{a&b\cr c&d}\right]\qquad{\rm and}\qquad \pi'=\pi
\left[\matrix{a&c\cr b&d}\right]$$
for any $a,b,c,d\in\{0,1\}$, forming a group isomorphic to
$\Z_2^4$. This  notation means
$$\pi\left[\matrix{a&b\cr c&d}\right](\la)\,\eqde \,J_v^{2a\,Q_v(\la)+2b\,Q_s(\la)}
J_s^{2c\,Q_v(\la)+2d\,Q_s(\la)}\la\ .\eqno(3.5)$$
When $k$ is odd, we will have six simple-current automorphisms: 
$$\eqalignno{
\pi&\,=\pi\left[\matrix{a&0\cr 0&d}\right]\qquad{\rm with}\qquad
\pi'=\pi\left[\matrix{a\,(d+1)&{dr\over 2}\cr {dr\over 2}&d}\right]&\cr
{\rm or}\quad \pi&\,=\pi\left[\matrix{{r\over 2}+1&b\cr
c&1}\right]\qquad{\rm with}\qquad
\pi'=\pi\left[\matrix{{r\over 2}+1+bc{r\over 2}&b+{r\over 2}\cr
{r\over 2}+1+bc+b&1}\right]\ ,&}$$
where  $a={r\over 2}$ or $d=0$, and where $b=1$ or $d=1$. The
corresponding permutation of $P_+$ is still given by (3.5). Again, for
these $r,k$, these are the values of $a,b,c,d$ for which (3.5) is
invertible. For $k$ odd, the group of simple-current automorphisms is isomorphic
to the symmetric group
${\frak S}_3$ when 4 divides $r$, and to $\Z_6$ when $r\equiv 2$ (mod 4).

For $k=2$ (so $\kappa=2r$), there are several Galois fusion-symmetries.
In particular, write $\la^i=\la^{2r-i}=\L_i$ for $1\le i\le
r-2$, and $\la^{r\pm 1}=\L_{r-1}+\L_r$. As with $B_{r,2}$,
$S_{00}^2={1\over 4\kappa}$
is rational so for any $m$ coprime to $2r$, we get a Galois
fusion-symmetry  $\pi\{
m\}$. It takes $\la^a$ to $\la^{ma}$, where the superscript is taken mod $2r$,
and will fix $J_v0$. Also,
$\pi\{m\}$ will send $J_s0$ to $J_s^m0$, as well as
stabilise the set $\{\L_r,\L_{r-1},J_v\L_r,J_v\L_{r-1}\}$. (In particular,
put $t=r$ when $r$ is even or when $m\equiv 1$ (mod 4), otherwise put
$t=r-1$; then for any $i,j$, $\pi\{m\}\,C_1^jJ_v^i\L_r$ is
$C_1^jJ_v^i\L_t$  or $C_1^jJ_v^{i+1}\L_t$, when the Jacobi symbol $({\kappa
\over m})$ is $\pm 1$, respectively.)\medskip

{\smcap Theorem 3.D.} {\it The fusion-symmetries of $D_r^{(1)}$ for 
$k\ne 2$
are all of the form $C_i\,\pi$, where $C_i$ is a conjugation,
and where $\pi$ is a simple-current automorphism. Similarly for $D_4^{(1)}$ at
$k=2$. Finally, when both $k=2$ and $r>4$, any 
fusion-symmetry $\pi$ can be written as $\pi=C_1^a\,\pi_v^b\,\pi\{m\}$
for $a,b\in\{0,1\}$ and any $m\in\Z_{2r}^{\times}$, $1\le m<r$.} 

\medskip $\pi_v$ here
refers to the simple-current automorphism $\pi[2]$ or $\pi[\matrix{1&0\cr0&0}]$,
 for $r$ odd/even. When $k=1$, ${\cal A}(D_{even,1})\cong {\frak S}_3$,
corresponding to any permutation of $\L_1,\L_{r-1},\L_r$, and $\A(D_{odd,1})
=\langle C_1\rangle\cong \Z_2$.
When $r>4$, $\A(D_{r,2})\cong (\Z^{\times}_{2r}/\{\pm 1\})\times \Z_2\times
\Z_2$ or $\Z_r^{\times}\times\Z_2$ for $r$ even/odd. 
$\A(D_{4,2})$ has 24 elements, and any element can be written
uniquely as $C_i\,\pi\left[\matrix{a&0\cr0&d}\right]$. 

\bigskip\noindent{\it 3.5. The algebra $E_6^{(1)}$}

\medskip
A weight $\la$ of
$P_+$ satisfies $k=\la_0+\la_1+2\la_2+3\la_3+2\la_4+\la_5+2\la_6$ and
$\kappa=k+12$. The charge-conjugation acts as $C\la
= (\la_0,\la_5,\la_4,\la_3,\la_2,\la_1,\la_6)$.
The order 3 simple-current $J$ is given by 
$J\la=(\la_5,\la_0,\la_6,\la_3,\la_2,\la_1,\la_4)$ with 
$Q(\la)=(-\la_1+\la_2-\la_4+\la_5)/3$.

The fusion products we need can be derived from [29] using (2.4):
$$\eqalignno{\L_1\stimes \L_1=&\,(\L_2)_2\splus(\L_5)_1\splus(2\L_1)_2&(3.6a)
\cr
\L_1\stimes\L_5=&\,(0)_1\splus(\L_6)_2\splus(\L_1+\L_5)_2&(3.6b)\cr
\L_1\stimes \L_2=&\,(\L_3)_3\splus(\L_6)_2\splus(\L_1+\L_2)_3\splus(\L_1+\L_5)_2
&(3.6c)\cr
\L_1\stimes (2\L_1)=&\,(3\L_1)_3\splus(\L_1+\L_2)_3\splus 
(\L_1+\L_5)_2&(3.6d)}$$ 
where the outer subscript on any summand denotes the smallest level where
that summand appears (it will also appear at all larger levels). So for example
$\L_1\stimes \L_1$ equals $\L_2\splus\L_5\splus (2\L_1)$ for any $k\ge 2$, but
equals $\L_5$ at $k=1$. A similar convention is used in (3.7) and elsewhere
 for higher fusion multiplicities (the number of
subscripts used will equal the numerical value of the fusion coefficient).\medskip

{\smcap Theorem 3.E6.} {\it The fusion-symmetries of $E_6^{(1)}$
are $C^i\,\pi[a]$, for any $i\in\{0,1\}$ and any $a\in\{0,1,2\}$ for which
$ak\not\equiv 1$ (mod $3$).}

\bigskip \noindent{\it 3.6. The algebra $E_7^{(1)}$}

\medskip
A weight $\la$ in $P_+$ satisfies $k=\la_0+2\la_1+3\la_2+4\la_3+3\la_4+2\la_5
+\la_6+2\la_7$, and $\kappa=k+18$. 
The charge-conjugation  is trivial, but there is a simple-current $J$ given
by $J\la=(\la_6,\la_5,\ldots,\la_1,\la_0,\la_7)$. It has $Q(\la)=
(\la_4+\la_6+\la_7)/2$.

The only fusion products we need can be obtained from [29] and (2.4):
$$\eqalignno{\L_6\stimes\L_6=&\,(0)_1\splus(\L_1)_2\splus(\L_5)_2\splus(
2\L_6)_2&\cr
\L_1\stimes \L_6=&\,(\L_6)_2\splus(\L_7)_2\splus(\L_1+\L_6)_3&\cr
\L_5\stimes\L_6=&\,(\L_4)_3\splus(\L_6)_2\splus(\L_7)_2\splus(\L_1+\L_6)_3
\splus(\L_5+\L_6)_3&\cr
\L_6\stimes(2\L_6)=&\,(\L_6)_2\splus(\L_1+\L_6)_3\splus(3\L_6)_3\splus(\L_5
+\L_6)_3&\cr
\L_4\stimes\L_6=&\,(\L_2)_3\splus(\L_3)_4\splus(\L_5)_3\splus(\L_1+\L_5)_4
\splus(\L_4+\L_6)_4\splus(\L_6+\L_7)_3&\cr
\L_6\stimes\L_7=&\,(\L_1)_2\splus(\L_2)_3\splus(\L_5)_2\splus (\L_6+\L_7)_3&\cr
\L_6\stimes(\L_5+\L_6)=&\,(\L_5)_3\splus(2\L_5)_4\splus(2\L_6)_3\splus(\L_6+\L_7)_3
\splus(\L_1+\L_5)_4&\cr&\,\splus(\L_4+\L_6)_4
\splus(\L_1+2\L_6)_4\splus(\L_5+2\L_6)_4&\cr}$$
 
At $k=3$ there is an order
3 Galois fusion-symmetry $\pi_3=\pi\{5\}$, which sends $J^i\L_1\mapsto J^i(2\L_6)
\mapsto J^i\L_2\mapsto J^i\L_1$ and fixes the other six weights.\medskip

{\smcap Theorem 3.E7.} {\it The only nontrivial fusion-symmetries for 
$E_7^{(1)}$ are $\pi[1]$  at even $k$, 
as well as $\pi_3$ and its inverse at $k=3$.}

\bigskip\noindent{\it 3.7. The algebra $E_8^{(1)}$}

\medskip
A weight $\la$ in $P_+$ satisfies $k=\la_0+2\la_1+3\la_2+4\la_3+5\la_4+6\la_5
+4\la_6+2\la_7+3\la_8$, and $\kappa=k+30$. 
The conjugations and  simple-currents are all trivial, except for  
an anomolous simple-current at $k=2$, sending $P_+=(0,\L_1,\L_7)$ to
$(\L_7,\L_1,0)$, which plays no role in this paper (except in Theorem 5.1).

The only fusion products we need can be derived from [28] and (2.4):
$$\eqalignno{\L_1\stimes\L_1=&\,(0)_2\splus(\L_1)_3\splus(\L_2)_3\splus
(\L_7)_2\splus(2\L_1)_4&(3.7a)\cr
\L_2\stimes\L_2=&\,(0)_3\splus(\L_1)_4\splus 2\sdot(\L_2)_{34}\splus 2\sdot
(\L_3)_{45}\splus (\L_4)_5\splus(\L_6)_4\splus 2\sdot
(\L_7)_{34}&(3.7b)\cr&\splus 2\sdot(\L_8)_{44}\splus 3\sdot(\L_1+\L_7)_{445}\splus
2\sdot(2\L_1)_{45}\splus(2\L_2)_6\!\splus (2\L_7)_4\!\splus
 2\sdot(\L_1+\L_2)_{55}
&\cr& \splus(\L_1+\L_3)_6\splus 2\sdot(\L_1+\L_8)_{55}\splus (\L_2+\L_7)_5
\splus (2\L_1+\L_7)_6\!\splus\!(3\L_1)_6&\cr
\L_7\stimes\L_7=&\,(0)_2\splus(\L_1)_3\splus(\L_2)_3\splus(\L_3)_4
\splus(\L_6)_4\splus(\L_7)_3\splus
(\L_8)_3&(3.7c)\cr &\,\splus(2\L_1)_4\splus(2\L_7)_4\splus(\L_1+\L_7)_4&\cr
(2\L_1)\!\stimes\!(2\L_1)=&\,(0)_4\splus(\L_1)_5\splus(\L_2)_5\splus(\L_3)_4
\splus(\L_7)_4\splus 2\sdot(2\L_1)_{46}\splus(2\L_2)_6&\cr &\splus
(2\L_7)_4
\splus 2\sdot(\L_1+\L_2)_{56}\splus(\L_1+\L_7)_5\splus(\L_2+\L_7)_5
\splus(3\L_1)_7 &\cr &
\splus(2\L_1+\L_2)_7\splus (2\L_1+\L_7)_6\splus(4\L_1)_8&(3.7d)\cr
\L_1\stimes\L_4=&\,(\L_3)_5\splus(\L_4)_6\!\splus(\L_5)_6\!\splus(\L_6)_5\!
\splus(\L_1+\L_3)_6\splus(\L_1+\L_4)_7\splus(\L_1+\L_6)_6&\cr&\splus\!(\L_1+\L_8)_5
\!\splus\!(\L_2+\L_7)_5\!\splus\!(\L_7+\L_8)_5\!\splus\!(\L_2+\L_8)_6\!\splus
\!(\L_3+\L_7)_6&(3.7e)\cr
\L_1\!\stimes\!(\L_1\!+\!\L_3)=&\,(\L_3)_6\splus(\L_4)_6\splus(\L_1+\L_2)_6\splus 
2\sdot(\L_1+\L_3)_{67}\splus(\L_1+\L_4)_7&\cr &\splus(\L_1+\L_6)_6
\splus(\L_1+\L_8)_6\splus(\L_2+\L_3)_7\splus(\L_2+\L_7)_6
\splus(2\L_2)_6&\cr&\splus(\L_2+\L_8)_6\splus(\L_3+\L_7)_6
\splus(2\L_1+\L_8)_7\splus
(2\L_1+\L_2)_7&\cr&\splus(2\L_1+\L_3)_8\splus(2\L_1+\L_7)_6\splus(\L_1+\L_2+\L_7)_7&(3.7f)\cr
\L_1\stimes(2\L_7)=&(\L_6)_4\!\splus\!(\L_1+\L_7)_4\!\splus\!(2\L_7)_5\!\splus
\!(\L_2+\L_7)_5\!\splus\!(\L_7+\L_8)_5\!\splus\!(\L_1+2\L_7)_6&(3.7g)\cr}$$

A fusion-symmetry at $k=4$, called $\pi_4$, was first found in [15]. It
interchanges
$\L_1\leftrightarrow\L_6$ and fixes the other eight weights in $P_+$. There also
is a fusion-symmetry, called $\pi_5$, at $k=5$ which
interchanges $\L_7\leftrightarrow 2\L_1$, $\L_8\leftrightarrow \L_1+\L_2$,
and $\L_6\leftrightarrow \L_2+\L_7$, and fixes the nine other weights.
The exceptional $\pi_5$  is closely related to the Galois permutation
$\la\mapsto\la^{(13)}$.

 \medskip
{\smcap Theorem 3.E8.} {\it The only nontrivial fusion-symmetries for $E_8^{(1)}$
are $\pi_{4}$ and $\pi_5$, occurring at $k=4$ and $5$ respectively.}

\bigskip \noindent{\it 3.8. The algebra $F_4^{(1)}$}

\medskip
A weight $\la$ in $P_+$ satisfies $k=\la_0+2\la_1+3\la_2+2\la_3+\la_4$,
and $\kappa=k+9$. Again, the conjugations and simple-currents are trivial. 

There are Galois fusion-symmetries at levels $k=3$ and 4. In particular, for $k=3$
we have the fusion-symmetry $\pi_3=\pi\{5\}$ which interchanges 
both $\L_2 \leftrightarrow\L_4$ and $\L_1\leftrightarrow 3\L_4$, and fixes
the other five weights in $P_+$. The exceptional $\pi_3$ was found independently
in [34,14]. For $k=4$ we get a
fusion-symmetry of order 4, which we will call $\pi_4$. It fixes 0, $\L_2+\L_4$,
$\L_3+\L_4$, and $2\L_4$, and permutes $\L_4\mapsto\L_1\mapsto 2\L_1\mapsto
4\L_4\mapsto\L_4$, $\L_2\mapsto 2\L_3\mapsto 3\L_4\mapsto \L_3\mapsto \L_2$,
and $\L_1+\L_3\mapsto \L_3+2\L_4\mapsto \L_1+\L_4\mapsto
\L_1+2\L_4\mapsto\L_1+\L_3$. Its square $\pi_4^2$ equals the 
fusion-symmetry $\pi\{5\}$. 

The only fusion products we need can be obtained from [29] and (2.4):
$$\eqalignno{\L_4\stimes\L_4=&\,(0)_1\splus(\L_1)_2\splus(\L_3)_2\splus
(\L_4)_1\splus(2\L_4)_2&\cr
\L_1\stimes\L_4=&\,(\L_3)_2\splus(\L_4)_2\splus(\L_1+\L_4)_3&\cr
\L_3\stimes \L_4=&\,(\L_1)_2\splus(\L_2)_3\splus(\L_3)_2\splus(\L_4)_2\splus
(\L_1+\L_4)_3\splus(\L_3+\L_4)_3\splus(2\L_4)_2&\cr
(2\L_4)\stimes\L_4=&\,(\L_3)_2\splus(\L_4)_2\splus(2\L_4)_2\splus(3\L_4)_3
\splus(\L_1+\L_4)_3\splus(\L_3+\L_4)_3&\cr}$$

{\smcap Theorem 3.F4.} {\it The only nontrivial fusion-symmetries of $F_4^{(1)}$
are $\pi_3$ at level $3$, and $\pi_4^i$ for $1\le i\le 3$,
which occur at level $4$.}

\bigskip 
\noindent{\it 3.9. The algebra $G_2^{(1)}$}\medskip

A weight $\la$ in $P_+$ satisfies $k=\la_0+2\la_1+\la_2$, and $\kappa=k+4$.
The conjugations and simple-currents are all trivial. 

Again there are nontrivial Galois fusion-symmetries.
At $k=3$, we have the order 3 fusion-symmetry $\pi_3=\pi\{4\}$ sending
$\L_1\mapsto 3\L_2\mapsto\L_2\mapsto \L_1$, and fixing the
remaining three weights. It was found in [14]. At $k=4$, we have $\pi_4=\pi\{5\}$
permuting both $\L_1\leftrightarrow 4\L_2$ and $2\L_1\leftrightarrow\L_2$,
and fixing the other five weights. It was found independently in [34,14],
and in \S5 we will see that it is closely related to the $\pi_3$ of $F_{4,3}$.

The only fusion products we will need can be obtained from [29] and (2.4):
$$\eqalignno{\L_2\stimes\L_2=&\,(0)_1\splus(\L_1)_2\splus(\L_2)_1\splus
(2\L_2)_2\ &\cr
\L_2\stimes\L_2\stimes\L_2=&\,(0)_1\splus 2\sdot(\L_1)_{22}\splus 4\sdot
(\L_2)_{1122}\splus 3\sdot(2\L_2)_{222}\splus 2\sdot(\L_1+\L_2)_{33}\splus
(3\L_2)_3   &\cr}$$

{\smcap Theorem 3.G2.} {\it The only nontrivial fusion-symmetries for $G_2^{(1)}$
are $(\pi_3)^{\pm 1}$ at $k=3$, and
$\pi_4$ at $k=4$.}

\bigskip\bigskip\centerline{{\smcap 4. The Arguments}} 
\bigskip

The fundamental reason the classification of fusion-symmetries for the affine algebras
is so accessible is (2.1b), which reduces the problem to studying Lie
group characters at elements of finite order. These values have been
studied by a number of people --- see e.g.\ [22,28] ---
 and the resulting combinatorics is often quite pretty.

Lemma 2.2 implies that a fusion-symmetry $\pi$ preserves q-dimensions:
$\D(\la)=\D(\pi\la)$ $\forall \la\in P_+$. In this subsection we use that
to find a weight $\L_\star$ for each algebra which must be essentially fixed by
$\pi$.

\vfill\eject\noindent{\it 4.1. q-dimensions}

\medskip  The most basic properties obeyed by the q-dimensions
$\D(\la)={S_{\la 0}\over S_{0 0}}$ 
are that $\D(\la)\ge 1$, and $\D(s\la)=\D(\la)$ for any $s\in\s$. Recall that
$\s$ is the symmetry group of the extended Dynkin diagram of $X_r^{(1)}$,
and that $s\in\s$ acts on $P_+$ by permuting the Dynkin labels.

The argument yielding Proposition 4.1 below relies heavily on the following observation.
Use (2.1c) to extend the domain of $\D$ from $P_+$ to 
the fundamental chamber $C_+$:
$$C_+\eqde \{\sum_{i=0}^rx_i\L_i\,|\,x_i\in\R,\ x_i>-1,\ \sum_{i=0}^r
x_ia_i^\vee=k\}\ .$$
Choose any $a,b\in C_+$. Then a straightforward calculation from (2.1c) gives
$${d\over dt}\D(ta+(1-t)b)=0\quad \Longrightarrow\quad {d^2\over dt^2}\D(ta+
(1-t)b)<0$$
for $0<t<1$. This means that for all $0<t<1$,
$$\D(ta+(1-t)b)>{\rm min}\{\D(a),\,\D(b)\}\ .\eqno(4.1)$$

\medskip
{\smcap Proposition 4.1} [17,18].\quad {\it For the following algebras $X_r^{(1)}$
and levels $k$, and choices of weight $\L_\star$, $\D(\la)=\D(\L_\star)$ implies
$\la\in\s\L_\star$:}

\item{(a)} {\it For $A_r^{(1)}$ any level $k$, where $\L_\star=\L_1$};

\item{(b)} {\it For $B_r^{(1)}$ any level $k\ne 2$, where $\L_\star=\L_1$};

\item{(c)} {\it For $C_r^{(1)}$ any level $k$ (except for $(r,k)=(2,3)$ or $(3,2)$),
where $\L_\star=\L_1$};

\item{(d)} {\it For $D_r^{(1)}$ any level $k\ne 2$, where $\L_\star=\L_1$};

\item{(e6)} {\it For $E_6^{(1)}$ any level $k$, where $\L_\star=\L_1$};

\item{(e7)} {\it For $E_7^{(1)}$ any level $k\ne 3$, where $\L_\star=\L_6$};

\item{(e8)} {\it For $E_8^{(1)}$ any level $k\ne 1,4$, where $\L_\star=\L_1$};

\item{(f4)} {\it For $F_4^{(1)}$ any level $k\ne 3,4$, where $\L_\star=\L_4$};

\item{(g2)} {\it For $G_2^{(1)}$ level any $k\ne 3,4$, where $\L_\star=\L_2$.}\medskip

The missing cases are: $B_{r,2}$ where $\D(\L_1)=\D(\L_2)=\cdots=
\D(\L_{r-1})=\D(2\L_r)$;

$D_{r,2}$ where $\D(\L_1)=\cdots=\D(\L_{r-2})$;

$C_{2,3}$ where $\D(\L_2)=\D(3\L_1)
=\D(\L_1)$, and its rank-level dual $C_{3,2}$;

$E_{7,3}$
where $\D(\L_1)=\D(\L_2)=\D(\L_6)$;

$E_{8,1}$ where $\L_1\not\in
P_+=\{0\}$, and  $E_{8,4}$ where $\D(\L_1)=\D(\L_6)$;

$F_{4,3}$ where $\D(\L_2)=\D(\L_4)$, and $F_{4,4}$ where
$\D(\L_1)=\D(2\L_1)=\D(4\L_4)=\D(\L_4)$;

$G_{2,3}$ where $\D(\L_1)=\D(\L_2)=\D(3\L_2)$, and $G_{2,4}$
where $\D(\L_2)=\D(2\L_1)$.\smallskip

The weight $\L_{\star}$ singled out by Proposition 4.1 (i.e.\ $\L_{\star}=
\L_1$ for $A_r^{(1)}$,
..., $\L_{\star}=\L_2$ for $G_2^{(1)}$)  is the nonzero weight 
with smallest Weyl dimension. What we find is that, for all but the smallest 
levels (see [18, Table 3]), $\L_{\star}$ will also have the smallest
q-dimension after the simple-currents.

The complete proof of Proposition 4.1 is given in [18], but to illustrate the
ideas we will sketch here the most interesting ($A_r^{(1)}$) and the
most difficult ($E_8^{(1)}$) cases.

Consider first $A_{r,k}$. By choosing $a-b=\L_i-\L_j$ in (4.1),
we get that either $\la=k\L_\ell$ for some $\ell$, in which case $\la$ is
a simple-current and (for $k\ne 1$) $\D(\la)<\D(\L_1)$, or $\D(\la)\ge
\D(\L_\ell)$ for some $\ell$, with equality iff $\la\in\s\L_\ell$. But then
rank-level duality $A_{r,k}\leftrightarrow A_{k-1,r+1}$ (defined as
for $C_{r,k}$, and which is exact for $A_{r,k}$ q-dimensions)
 and (4.1) with $a-b=\widetilde{\L_0}-\widetilde{\L_1}$
give us $\D(\L_\ell)=\widetilde{\D}(\ell\widetilde{\L_1})\ge \widetilde{\D}
(\widetilde{\L_1})=\D(\L_1)$, with equality iff $\ell=1$ or $r$. Combining
these results yields Proposition 4.1(a).

For $E_{8,k}$, run through each $a-b=a_j^\vee\L_i-a_i^\vee\L_j$ to reduce
the proof to comparing $\D(\L_1)$ with $\D({k\over a_i^\vee}\L_i)$ for
$i\ne 0$, or $\D(\L_i)$ for $i\ne 0,1$ (the argument in [18] unnecessarily
complicated things by restricting to integral weights). Standard arguments
(see [18] for details) quickly show that the q-dimension
$\D({k\over a_i^\vee}\L_i)$ monotonically increases with $k$ to $\infty$,
while $\D(\L_i)$ monotonically increases with $k$ to the Weyl dimension of
$\L_i$. The proof of Proposition 4.1(e8) then reduces
to a short computation.

\bigskip\noindent{{\it 4.2. The $A$-series argument}}

\medskip
Recall that $\r=r+1$.
Proposition 4.1(a) tells us that $\pi\L_1=C^{a}J^{b}\L_1$, for
 some $a,b$. Hitting $\pi$ with $C^a$, we can assume without loss of
  generality that $a=0$.
Write $\pi(J0)=J^c0$; then $\pi$ can be a permutation of $P_+$ only if
$c$ is coprime to $\r$.

If $k=1$ then $P_+=\{0,J0,\ldots,J^r0\}$ so $\pi=\pi[c-1]$. Thus we can assume
$k\ge 2$.

Useful is the coefficient of $\la$ in the tensor product $\L_1\otimes\cdots
\otimes\L_1$ ($\ell$ times): it is 0 unless $t(\la)=\ell$, in which case
the coefficient is ${\ell !\over h(\la)}$ (to get this, compare (3.1)
above with [27, p.114]) --- we equate here the fundamental
weights $\L_{\r}$ and $\L_0$, so e.g.\ `${k\over\r}\L_{\r}$' equals
`0' when $\r$ divides $k$. Here, $h(\la)=\prod h(x)$ is the
product of the hook-lengths of the Young diagram corresponding to $\la$.
Equation (2.4) tells us that as long as $t(\la)=\ell\le k$, the
number ${\ell !\over h(\la)}$ will also be the coefficient of $N_\la$ in
the fusion
power $(N_{\L_1})^\ell$. Note that $J0=k\L_1$ is the only simple-current
appearing in the fusion product $\L_1\stimes\cdots\stimes\L_1$ ($k$ times).
Thus the only nontrivial simple-current appearing in the fusion $\pi\L_1\stimes\cdots
\stimes\pi\L_1$ will be $J^{bk}J0$ (0 will appear iff $\r$ divides
$k$).  Hence $bk+1\equiv c$ (mod $\r$) must
be coprime to $\r$. This is precisely the condition needed for $\pi[b]$ to
be a simple-current automorphism.

In other words, it suffices to consider $\pi\L_1=\L_1$ and hence $\pi[J0]=
J0$. We are done if $r=1$, so assume $r\ge 2$. 
 From the $\L_1\stimes\L_1$ fusion, we get
that $\pi\L_2\in\{\L_2,2\L_1\}$. Note that $k\L_1$ occurs (with multiplicity 1)
in the tensor and fusion product of $2\L_1$ with $k-2$ $\L_1$'s, but that
it doesn't in the tensor (hence fusion) product of $\L_2$ with $k-2$ $\L_1$'s
(recall that $k\L_1\succ (k-2)\L_1+\L_2$ in the usual partial order on
weights). Since $\L_2\stimes\L_1\stimes\cdots\stimes\L_1$ does not contain $J0$,
$(\pi\L_2)\stimes(\pi \L_1)\stimes\cdots\stimes(\pi\L_1)$ should also
avoid $\pi(J0)=J0$, and thus $\pi\L_2$ cannot equal $2\L_1$.

Thus we know $\pi\L_2=\L_2$. The remaining $\pi\L_\ell=\L_\ell$ follow
quickly from induction: if $\pi\L_\ell=\L_\ell$ for some $2\le \ell<r$, then
the fusion $\L_1\stimes\L_\ell$ tells us
$\pi\L_{\ell+1}\in\{\L_{\ell+1},\L_1+\L_\ell\}$. But $h(\L_1+\L_\ell)=
(\ell+1)!/\ell$ and $h(\L_{\ell+1})=(\ell+1)!$, so $\pi\L_{\ell+1}=\L_{\ell+1}$.
Thus $\pi$ fixes
all fundamental weights, and since these comprise a fusion-generator (see
the discussion at the end of \S2.2) we know that $\pi$ must fix everything
in $P_+$.

\vfill\eject\noindent{{\it 4.3. The $B$-series argument}}

\medskip
$k=1$ is easy: $P_+=\{0,J0,\L_r\}$ and $\pi=id.$ is automatic. $k=2$
will be done later in this subsection. Assume now that $k\ge 3$.

 From Proposition 4.1(b) we can
write $\pi\L_1=J^{a}\L_1$ and $\pi'\L_1=J^{a'}\L_1$. We know $\pi J0=J0$,
 so (2.7b) says $\pi$ must
take spinors to spinors, and nonspinors to nonspinors. Then we will have
$\chi_{\L_1}[\psi]=(-1)^{a'}\chi_{\L_1}[\pi\psi]$ for any spinor $\psi$.
 Now if $a'=1$, then $\pi$ will take the
spinors which maximize $\chi_{\L_1}$, to those which minimize it. Both these maxima
and minima can be easily found from (3.2). Thus we get that
$\pi(\s\L_r)$ equals $k\L_r$ (when $k$ odd) or $\s((k-1)\L_r)$ (when $k$ even).
But the sets $\s\L_r$ and $k\L_r$ have different cardinalities ($k\L_r$ is
a $J$-fixed-point), and so can't get mapped to each other. Also, the fusions
$\L_1\stimes\L_r=\L_r\splus(\L_1+\L_r)$ and $J^a\L_1\stimes (J^i(k-1)\L_r)
=(J^{a+i}(k-1)\L_r)\splus(J^{a+i+1}(k-1)\L_r)\splus J^{a+i+1}(\L_{r-1}+(k-3)
\L_r)$ have different numbers of weights on their right sides, so also
$\pi\L_r\not\in\s(k-1)\L_r$.

Thus $a'=0$ and $\pi\L_r=J^b\L_r$ for some $b$. Similarly, $a=0$.
Hitting $\pi$ with $\pi[1]^b$, we  may assume that $\pi$ fixes
 $\L_r$.

Now assume $\pi$ fixes $\L_\ell$, for $1\le\ell<r-1$. Then the fusion 
$\L_1\stimes\L_\ell$ says that $\pi\L_{\ell+1}$ equals $\L_{\ell+1}$ or
$\L_1+\L_\ell$. But from (3.2) we find 
$$\eqalign{\chi_{\L_1}[\L_{\ell+1}]-\chi_{\L_1}[\L_1+\L_\ell]=&\,
2\,\{\cos(\pi{2r-2\ell+1
\over\kappa})-\cos(\pi{2r-2\ell-1\over\kappa})+\cos(\pi{2r+1\over\kappa})
\cr-&\cos(\pi{ 2r+3\over\kappa})\}=4\cos(\pi{2r-\ell+1\over\kappa})\,
\{\cos(2\pi{\ell\over\kappa})-\cos(2\pi{\ell+1\over\kappa})\}>0}$$
Hence $\pi$ will
 fix $\L_{\ell+1}$ if it fixes $\L_\ell$, concluding the
argument.             \smallskip

Now consider the more interesting case: $k=2$. Then $\kappa=2r+1$; recall
the weights in $P_+(B_{r,2})$ are the simple-currents $0,J0$, the $J$-fixed-points
$\ga^1,\ldots,\ga^r$ (notation defined in \S3.2), and the spinors $\L_r,J\L_r$.
Because $\pi(J0)=\pi'(J0)=J0$, we know both $\pi$ and $\pi'$ must take $J$-fixed-points
to $J$-fixed-points, i.e.\ $\pi
\L_1=\ga^{m}$ and $\pi'\L_1=\ga^{m'}$ for some $1\le m,m'\le r$. It is easy to compute [25]
$${S_{\ga^a\ga^b}\over S_{0\ga^b}}=2\cos(2\pi {ab\over \kappa})\ .\eqno(4.2)$$ 
 From this we see $m\,m'\equiv \pm 1$ (mod $\kappa$), so $m$ is coprime to $\kappa$.
Hitting it with the Galois fusion-symmetry $\pi\{m'\}$, we see that we may
assume $\pi\L_1=\pi'\L_1=\L_1$.

Now use (4.2) to get $\pi\ga^i=\ga^i$ for all $i$. Then $\pi$ equals
the identity or $\pi[1]$, depending on what $\pi$ does to $\L_r$.

\bigskip\noindent{{\it 4.4.\ The $C$-series argument}}

\medskip By rank-level duality,
we may take $r\le k$. For now assume $(r,k)\ne(2,3)$. Then we know $\pi
\L_1=J^{a}\L_1$ and $\pi\L_1=J^{a'}\L_1$ for some $a,a'$. Since $\pi J0=\pi'
J0=J0$, (2.7b) says $a=a'=0$ if $kr$ is odd. Since
$\chi_{\L_1}[\L_1]>0$ (using (3.3)), $S_{\L_1\L_1}=S_{J^a\L_1,J^{a'}\L_1}$
implies that $a=a'$ also holds when $kr$ is even, and hence we may
assume (hitting with $\pi[1]^a$) that also $a=a'=0$ holds for $kr$
even.  From the fusion $\L_1
\stimes\L_\ell$ we get $\pi\L_{\ell+1}\in\{\L_{\ell+1},\L_1+\L_\ell\}$
if $\pi\L_\ell=\L_\ell$; for $r<k$ conclude the argument with the calculation
$$\chi_{\L_1}[\L_{\ell+1}]-\chi_{\L_1}[\L_1+\L_\ell]=4\cos(\pi\,{2r+2-\ell\over
2\kappa})\,\{\cos(\pi\,{\ell\over 2\kappa})-\cos(\pi\,{\ell+2\over
2\kappa})\}>0$$
as in \S4.3. When $r=k$, that 
inequality only holds for $\ell>1$, but we can force $\pi\L_2=\L_2$ by
hitting $\pi$ if necessary with $\pi_{{\rm rld}}$.

The remaining case $C_{2,3}$ follows because $\pi'J0=J0$: by (2.7b)
$\pi\L_1\not\in\s\L_2$, and by (2.7a) $\pi\L_1\ne 3\L_1$ ($3\L_1$ is
a $J$-fixed-point).

\bigskip\noindent{{\it 4.5.\ The $D$-series argument}}

\medskip $k=1$ is trivial,
and $k=2$ will be considered shortly. For $k>2$, Proposition 4.1 tells us
that $\pi\L_1=J_v^aJ_s^b\L_1$ and $\pi'\L_1=J_v^{a'}J_s^{b'}\L_1$, for
$a,a',b,b'\in\{0,1\}$. Immediate from (3.4) is that $\chi_{\L_1}
[\L_1]>0$ and that $\chi_{\L_1}[\psi]$, for a spinor $\psi$, takes its
maximum at $C^iJ_v^j\L_r$. Our first step is to force
$\pi\L_1=\pi'\L_1=\L_1$. Unfortunately this requires a case analysis.

Consider first even $r\ne 4$, and even $k>2$. Now, $0\ne S_{\L_1\L_1}=
S_{\pi\L_1,\pi'\L_1}$ forces $b=b'$; hence
hitting with the simple-current automorphism $\pi\left[\matrix{0&a\cr a'&b\cr}
\right]$, we may assume $\pi\L_1=\pi'\L_1=\L_1$.

Next consider even $r\ne 4$ and odd $k>2$. Either of $\pi\L_1=J_v\L_1$ or
$\pi'\L_1=J_v\L_1$
is impossible, by comparing $S_{\L_1,J_s0}$ and $S_{J_v\L_1,J0}$ for any
simple-current $J$. For any of the three remaining choices of $J_v^aJ_s^b\L_1$,
we can find a simple-current automorphism of the form $\pi\left[\matrix{*&a\cr
*&b\cr}\right]$; hitting its inverse onto $\pi$ allows us to take $a=b=0$.
Again $0\ne S_{\L_1\L_1}$ forces $b'=0$, and now $a'=1$ is forbidden.
Thus again $\pi\L_1=\pi'\L_1=\L_1$.

As usual, $r=4$ is complicated by triality. We can force $\pi\L_1=\L_1$.
That we can also take $\pi'\L_1=\L_1$, follows from the inequality
$\chi_{\L_1}[\L_1]>\chi_{\L_1}[\L_3]=\chi_{\L_1}[\L_4]>0$, valid
for $k\ge 3$. Establishing that inequality from (3.4) is equivalent to
showing
$$1+\cos(x)+\cos(2x)+\cos(4x)>\cos(x/2)+\cos(3x/2)+\cos(5x/2)+\cos(7x/2)$$
for $0<x\le 2\pi/9$, which can be shown e.g.\ using Taylor series.

For odd $r$, the charge-conjugation $C$ equals $C_1$. Since it must commute with
 $\pi$, i.e.\ that $C_1\pi\L_1=J_v^{a+b}J_s^b\L_1$ must equal $\pi C_1\L_1
 =J_v^aJ_s^b\L_1$, we get that $b=0$. Similarly $b'=0$. When $k$ is odd,
eliminate $a=1$ and $a'=1$ by comparing $S_{\L_1,J_s0}$ and $S_{J_v\L_1,J0}$
as before. The hardest case is $k$ even. We can force $\pi\L_1=\L_1$ by hitting
with $\pi[a]$. Suppose for contradiction that $\pi'\L_1=J_v\L_1$. We
 know $\pi'(J_v0)=J_v0$ (compare $S_{\L_1,J_v0}$ and $S_{\L_1,
J0}$), so by (2.7b) $\pi\L_r$ must be  a spinor. $\chi_{\L_1}[\L_r]=\chi_{J_v\L_1}[
\pi\L_r]$ requires $\pi\L_r=C_1^iJ_v^jJ_s\L_r$. From the $\L_1\stimes\L_r$
fusion we get $\pi\L_{r-1}=C_1^iJ_v^jJ_s\L_{r-1}$, but $C\pi=\pi C$ says
that $\pi\L_{r-1}=C_1^iJ_v^{j+1}J_s\L_{r-1}$ --- a contradiction.

Thus in all cases we have $\pi\L_1=\pi'\L_1=\L_1$.
We know $\pi'(J_v0)=J_v0$ (compare $S_{\L_1,J_v0}$ and $S_{\L_1,J0}$),
so  $\pi\L_r$ is a spinor
and in fact must equal $\pi\L_r=C_1^iJ_v^j\L_r$. Hitting with
$(C_1^i\pi_v^j)^{-1}$, 
 we can require $\pi\L_r=\L_r$. That $\pi\L_{r-1}$ must now
equal $\L_{r-1}$ follows from the $\L_1\stimes\L_r$ fusion.

Next, note that we know from $\L_1\stimes\L_1$ that $\pi\L_2$ is
$\L_2$ or $2\L_1$. As in \S4.2, the fusion
$(2\L_1)\stimes\L_1\stimes\cdots \stimes\L_1$ ($k-2$ times) contains
the simple-current $J_v0$, but
$\L_2\stimes\L_1\stimes\cdots\stimes\L_1$ ($k-2$ times) doesn't, and
thus $\pi\L_2=\L_2$.

Assume $\pi\L_\ell=\L_\ell$. Using the fusions $\L_1\stimes\L_\ell$ (for $1<\ell<r-2$), and noting that 
$$\chi_{\L_1}[\L_{\ell+1}]-\chi_{\L_1}[\L_1
+\L_\ell]=4\cos(\pi\,{2r-\ell\over\kappa})\,\{\cos(\pi\,{\ell\over\kappa})-
\cos(\pi\,{\ell+2\over\kappa})\}$$
equals 0 only when $\ell=r+1-k/2$, we see that
$\pi\L_{\ell+1}=\L_{\ell+1}$ except possibly for $\ell=r+1-k/2$ (hence
$2r-2\ge k\ge 4$).
For that $\ell$, use q-dimensions:
$${\D(\L_1+\L_\ell)\over \D(\L_{\ell+1})}={\sin(2\pi\,(k-2)
/\kappa)\over \sin(2\pi /\kappa)}>1\ ,$$
which is valid for these $k$. 
So we also know $\pi\L_{i}=\L_{i}$ for all $i\le r-2$, and
we are done.

\smallskip All that remains is $D_{r,2}$. Recall the $\la^i$ defined in \S3.4.
 Note that $\D(\L_r)=\sqrt{r}$, $\D(\la^a)=2$, and $S_{\la^a\la^b}
/S_{0\la^b}=2\cos(\pi ab/r)$. For $r\ne 4$, the q-dimensions force $\pi
\L_1=\la^{m}$ and $\pi'\L_1=\la^{m'}$, and $S_{\L_1\L_1}=S_{\la^m\la^{m'}}$
says $mm'\equiv \pm 1$ (mod $2r$).
So without loss of generality we may take $m=m'=1$. The rest of the argument
is easy.

For $D_{4,2}$, we can force $\pi\L_1=\L_1$, and then eliminate
$\pi'\L_1=\L_{r-1}$ or $\L_r$ by $S_{\L_1\L_1}\ne 0=S_{\L_1\L_r}=S_{\L_1
\L_{r-1}}$. The rest of the argument is as before.

\bigskip
\noindent{{\it 4.6. The arguments for the exceptional algebras}}

\medskip
The exceptional algebras follow quickly from the fusions (and Dynkin
diagram symmetries) given in \S\S3.5-3.9.

For example, consider $E_6^{(1)}$ for $k\ge 2$. Proposition 4.1 tells us
$\pi\L_1=C^aJ^b\L_1$ for some $a,b$, and we know $\pi'J0=J^c0$ for $c=\pm 1$.
Hence from (2.7b) we get $kb\not\equiv -1$ (mod
3). Hitting $\pi$ with $\pi[-b]^{-1}
C^a$, we need consider only $\pi\L_1=\L_1$. It is now immediate that
 $\pi\L_5=\L_5$, by
commuting $\pi$
with $C$. From (3.6a) we get that $\pi$ must permute $\L_2$ and $2\L_1$. Compare
(3.6c) with (3.6d): since for any $k\ge 2$ they have different numbers of summands,
we find in fact that $\pi$ will fix both $\L_2$ (hence $\L_4$) and $2\L_1$. 
 From (3.6b) we get that $\pi$ permutes $\L_6$ and $\L_1+\L_5$, and so (3.6d)
now tells us $\pi\L_6=\L_6$. Finally, (3.6c) implies (for $k\ge 3$)
$\pi\L_3=\L_3$ (since $C\pi=\pi C$), and we are done for $k\ge
3$. Since $\{\L_1,\L_2,\L_4,\L_5,\L_6\}$ is a fusion-generator for
$k=2$ (see \S2.2), we are also done for $k=2$. 

For $E_8^{(1)}$ when $k\ge 7$, (3.7a) tells us that $\L_2,\L_7,2\L_1$
are permuted. For those $k$, the highest multiplicities in (3.7b)--(3.7d)
 are 3, 1, 2, respectively, so $\L_2,\L_7,2\L_1$
must all be fixed. The fusion product (3.7c) also tells us
that $\L_3,\L_6,\L_8,\L_1+\L_7,2\L_7$ are permuted; (3.7d)
then says that the sets $\{\L_6,\L_8\}$, $\{\L_3,\L_1+\L_7,2\L_7\}$, and $\{2\L_2,
\L_2+\L_7,3\L_1,2\L_1+\L_2,2\L_1+\L_7,4\L_1\}$ are stabilised.
Now (3.7b) implies $\L_3,\L_6,\L_8,2\L_7$ are all fixed, while
the set $\{\L_4,\L_1+\L_3\}$ is stabilised. Comparing (3.7e) and (3.7f),
we get that $\L_4$ is fixed and $\L_5,\L_7+\L_8$
are permuted. Finally, (3.7g) shows $\L_5$ also is fixed.
To do $E_8^{(1)}$ when $k\le 6$, knowing q-dimensions really simplifies
things.

\vfill\eject\centerline{{\bf 5. Affine fusion ring isomorphisms}}\bigskip

We conclude the paper with the determination of all isomorphisms among the
affine fusion rings ${\cal R}(X_{r,k})$. Recall Definition 2.1 and the
discussion in \S2.2.

\medskip{\smcap Theorem 5.1.} {\it The complete list of fusion ring isomorphisms
 ${\cal R}(X_{r,k})\cong{\cal R}(Y_{s,m})$ when $X_{r,k}\ne Y_{s,m}$ (where
 $X_r,Y_s$ are simple) is: 

\noindent rank-level duality ${\cal R}(C_{r,k})\cong{\cal R}(C_{k,r})$ for all
 $r,k$, as well as ${\cal R}(A_{1,k})\cong{\cal R}(C_{k,1})$;

 \noindent ${\cal R}(B_{r,1})\cong{\cal R}(A_{1,2})
 \cong{\cal R}(C_{2,1})\cong{\cal R}(E_{8,2})$ for all $r\ge 3$;

 \noindent ${\cal R}(A_{3,1})\cong{\cal R}(D_{odd,1})$;

 \noindent ${\cal R}(D_{r,1})
 \cong{\cal R}(D_{s,1})$ whenever $r\equiv s$ (mod $2$);

\noindent ${\cal R}(A_{2,1})\cong
 {\cal R}(E_{6,1})$ and ${\cal R}(A_{1,1})\cong{\cal R}(E_{7,1})$;

\noindent ${\cal R}
 (F_{4,1})\cong{\cal R}(G_{2,1})$, ${\cal R}(F_{4,2})\cong
{\cal R} (E_{8,3})$, and ${\cal R} (F_{4,3})\cong{\cal R}(G_{2,4})$.}

\medskip The isomorphism ${\cal R}(A_{1,k})\cong{\cal R}(C_{k,1})$
takes $a\L_1$ to $\widetilde{\L}_a$. The isomorphism 
${\cal R}(F_{4,2})\cong{\cal R}(E_{8,3})$ was found in [14]; 
it relates $\L_1\leftrightarrow\widetilde{\L}_8$,
$2\L_4\leftrightarrow \widetilde{\L}_2$, $\L_3\leftrightarrow\widetilde{\L}_1$,
$\L_4\leftrightarrow\widetilde{\L}_7$. The isomorphism 
${\cal R} (F_{4,3})\cong{\cal R}(G_{2,4})$ was found in [34,14]; a
correspondence which works is $\L_4\leftrightarrow\widetilde{\L}_1$,
$\L_1\leftrightarrow 2\widetilde{\L}_1$, $\L_3\leftrightarrow
3\widetilde{\L}_2$, $2\L_4\leftrightarrow 2\widetilde{\L}_2$, 
$\L_1+\L_4\leftrightarrow\widetilde{\L}_1+2\widetilde{\L}_2$,
$\L_2\leftrightarrow 4\widetilde{\L}_2$, $3\L_4\leftrightarrow\widetilde{\L}_2$,
and $\L_3+\L_4\leftrightarrow\widetilde{\L}_1+\widetilde{\L}_2$.

 We will sketch the proof here. 
The idea is to compare invariants for the various fusion rings, case by case.
For example, suppose ${\cal R}(A_{r,k})$ and ${\cal R}(A_{s,m})$ are
isomorphic. Then their simple-current groups $\Z_{r+1}$ and $\Z_{s+1}$
must be isomorphic (since simple-currents must get mapped to simple-currents), so
$r=s$. Now compare the numbers $\|P_+\|$ of highest-weights: $({r+k\atop r})=
({r+m\atop r})$, which forces $m=k$.

It is also quite useful here to know those weights with second smallest
q-dimension. This is a by-product of the proof of Proposition 4.1, and
the complete answer is given in [18, Table 3]. Here we will simply
state that those weights in $P_+^k(X_r^{(1)})$ with second smallest
q-dimension are precisely the orbit $\s\L_\star$, except for: $A_{r,1}$;
$B_{r,k}$ for $k\le 3$; $C_{2,2},C_{2,3},C_{3,2}$; $D_{r,k}$ for $k\le 2$;
$E_{6,k}$ for $k\le 2$; and $E_{7,k},E_{8,k},F_{4,k},G_{2,k}$ for $k\le
4$.

$C_{r,k}$ and $B_{s,m}$ both have two simple-currents, but their fusion rings
can't be isomorphic (generically) because the orbit
$J^i\L_1$ has the second smallest q-dimension for both algebras at
generic rank/level,
but the numbers $Q_j(J^i\L_1)$ for the two algebras are different.

Another useful invariant involves the set of integers $\ell$ coprime
to $\kappa N$ for which
$0^{(\ell)}$ is a simple-current. For the classical algebras this is easy to find,
using (2.1c):
Up to a sign, the q-dimension of $0^{(\ell)}$ ($\ell$ coprime to $2\kappa$)
for the algebras $B_r^{(1)},C_r^{(1)},D_r^{(1)}$ is, respectively,
$$\eqalign{
\prod_{a=0}^{r-1}&\;{\sin(\pi\ell\,(2a+1)/2\kappa)\over \sin(\pi\,(2a+1)/2\kappa)}\,
\prod_{b=1}^{2r-2}{\sin(\pi\ell b/\kappa)^{[{2r-b\over 2}]}\over\sin(\pi
  b/\kappa)^{[{2r-b\over 2}]}}\ ,\cr
\prod_{a=1}^{r-1}&\;{\sin(\pi\ell a/\kappa)^{r-a}\,\sin(\pi\ell\,(2a-1)/2\kappa)^{r-a}
\over \sin(\pi a/\kappa)^{r-a}\,\sin(\pi\,(2a-1)/2\kappa)^{r-a}}\,
\prod_{b=r}^{2r-1}{\sin(\pi\ell b/2\kappa)\over \sin(\pi b/2\kappa)}\ ,\cr
\prod_{a=1}^{r-1}&\;{\sin(\pi\ell a/\kappa)^{[{2r-a+1\over 2}]}\over \sin(\pi
  a/\kappa)^{[{2r-a+1\over 2}]}}\,\prod_{b=r}^{2r-3}{\sin(\pi\ell
    b/\kappa)^{[{2r-b-1\over 2}]} \over \sin(\pi
b/\kappa)^{[{2r-b-1\over 2}]}} \ ,\cr}$$
where $[x]$ here denotes the greatest integer not more than $x$.
The absolute value of each of these is quickly seen to be greater than 1 unless $\ell
\equiv \pm 1$ (mod $2\kappa$), except for the orthogonal algebras
when $k\le 2$. An isomorphism ${\cal R}(X_{r,k})
\cong{\cal R}(X_{r',k'})$ would require then that whenever
$\ell\equiv\pm 1$ (mod $2\kappa$) is coprime to $\kappa'$, it must
also satisfy $\ell\equiv\pm 1$ (mod $2\kappa'$), and conversely. This
forces $\kappa= \kappa'$, for $X=B$ or $D$ and $k>2$, or $X=C$ and any
$k$.

If ${\cal R}(C_{r,k})\cong {\cal R}(C_{s,m})$, then that Galois argument
implies $r+k+1=
s+m+1$, so compare numbers of highest-weights: $({r+k\atop r})=({r+k\atop s})$.

A similar argument works for the orthogonal algebras. For instance suppose
 ${\cal R}(B_{r,k})\cong{\cal R}(B_{s,m})$ but $B_{r,k}\ne B_{s,m}$, and that
 $k,m>2$. Then Galois
 implies  $2r+k=2s+m$. Comparing the value of ${\cal D}(\L_1)$ (the second
 smallest q-dimension when $k>3$), using (3.2) with $\la=0$, tells us that
 $2s+1=k,2r+1=m$. Now count the number of fixed-points of $J$ in both
 cases: $({\kappa/2-1\atop r-1})=({\kappa/2-1\atop s-1})$, i.e.\ $s-1=(k-1)/2$,
 a contradiction.

For comparing classical algebras with exceptional algebras, a useful device
is to count the number of weights appearing in the fusion
$\L_\star\stimes\L_\star$ (when $\L_\star$ has second smallest
q-dimension).  For example, for $A_{1,k}$ ($k>1$), $C_{r,k}$ ($k>1$,
except for $C_{2,2},C_{2,3},C_{3,2}$), and
$E_{7,k}$ ($k>4$), we learned in \S3 that this number is 2, 3, 4 respectively,
so none of these can be isomorphic.

For the orthogonal algebras at level 2, useful is the number
of weights with second smallest q-dimension (respectively $r$ and $r-1$
for $B_{r,2}$ and $D_{r,2}$, except for $D_{4,2}$).

For the exceptional algebras, comparing ${\cal D}(\L_\star)$
and the number of highest-weights is effective. Recall that both $\|P_+\|$ and
${\cal D}(\L_\star)$
for a fixed algebra monotonically increase with $k$ to (respectively)
$\infty$ and the Weyl dimension
of $\L_\star$, which is 7, 26, and 248 for $G_2,F_4,E_8$ respectively.
For $E_{8,k}$, ${\cal D}(\L_1)$ exceeds 7 for $k\ge 5$, and exceeds 26 for
$k\ge 11$, while $F_{4,k}$ exceeds 7 for $k\ge 4$. The number of highest-weights
of $E_{8,4},E_{8,10}$, and $F_{4,3}$ are 10, 135, and 9, so only a small
number of possibilities need be considered.

\bigskip \centerline{{\smcap References}} \medskip

\item{1.} {\smcap R.\ E.\ Behrend, P.\ A.\ Pearce, V.\ B.\ Petkova}
and {\smcap J.-B.\ Zuber}, On the classification of bulk and boundary
conformal field theories, {\it Phys.\ Lett.} {\bf B444} (1998),
163--166;

\item{} {\smcap J.\ Fuchs} and {\smcap C.\ Schweigert}, Branes: From
free fields to general backgrounds, {\it Nucl.\ Phys.} {\bf B530} (1998),
99--136.

\item{2.} {\smcap D.\ Bernard}, String characters from Kac--Moody automorphisms,
{\it Nucl.\ Phys.} {\bf B288} (1987), 628--648.

\item{3.} {\smcap J.\ B\"ockenhauer} and {\smcap D.\ E.\ Evans},
Modular invariants from subfactors: Type I coupling matrices and
intermediate subfactors, preprint math.OA/9911239, 1999.

\item{4.} {\smcap A.\ Coste} and {\smcap T.\ Gannon}, Remarks on
Galois in rational conformal field theories, {\it Phys.\ Lett.} {\bf B323} (1994),
316--321.

\item{5.} {\smcap A.\ Coste, T.\ Gannon} and {\smcap P.\ Ruelle},
Finite group modular data, preprint hep-th/0001158, 2000.

\item{6.} {\smcap Ph.\ Di Francesco, P.\ Mathieu} and {\smcap D.\
S\'en\'echal}, ``Conformal Field Theory'', Springer-Verlag, New York, 1997.

\item{7.} {\smcap G.\ Faltings}, A proof for the Verlinde formula,
{\it J.\ Alg.\ Geom.} {\bf 3} (1994), 347--374.

\item{8.} {\smcap M.\ Finkelberg}, An equivalence of fusion categories,
{\it Geom.\ Func.\ Anal.} {\bf 6} (1996), 249--267.

\item{9.} {\smcap I.\ B.\ Frenkel}, Representations of affine Lie algebras,
Hecke modular forms and Korteg-de Vries type equations, in: ``Lie
algebras and related topics'', Lecture
Notes in Math, Vol.\ 933, Springer-Verlag, New York, 1982.

\item{10.} {\smcap I.\ B.\ Frenkel, Y.-Z.\ Huang} and {\smcap J.\
Lepowsky}, ``On axiomatic approaches to vertex operator algebras and modules'',
{\it Memoirs Amer.\ Math.\ Soc.} {\bf 104} (1993).

\item{11.} {\smcap J.\ Fr\"ohlich} and {\smcap T.\ Kerler}, ``Quantum
groups,  quantum categories and quantum field theory'', Lecture Notes
in Mathematics, Vol.\ 1542, Springer-Verlag, Berlin, 1993.

\item{12.} {\smcap J.\ Fuchs}, Simple WZW currents, {\it Commun.\ Math.\ Phys.}
{\bf 136} (1991), 345--356.

\item{13.} {\smcap J.\ Fuchs, B.\ Gato-Rivera, A.\ N.\ Schellekens}
and {\smcap C.\ Schweigert}, Modular invariants and fusion automorphisms
from Galois theory, {\it Phys.\ Lett.} {\bf B334} (1994), 113--120.

\item{14.} {\smcap J.\ Fuchs} and {\smcap P.\ van Driel}, WZW fusion
rules, quantum groups, and the modular matrix $S$, {\it Nucl.\ Phys.} {\bf B346}
(1990), 632--648.

\item{15.} {\smcap J.\ Fuchs} and {\smcap P.\ van Driel}, Fusion rule
engineering, {\it Lett.\ Math.\ Phys.} {\bf 23} (1991), 11--18.

\item{16.} {\smcap T.\ Gannon}, WZW commutants, lattices, and level-one
partition functions, {\it Nucl.\ Phys.} {\bf B396} (1993), 708--736;

\item{} {\smcap P.\ Ruelle, E.\ Thiran} and {\smcap J.\ Weyers}, Implications
of an arithmetical symmetry of the commutant for modular invariants,
{\it Nucl.\ Phys.} {\bf B402} (1993), 693--708.

\item{17.} {\smcap T.\ Gannon}, Symmetries of the Kac-Peterson modular 
matrices of affine algebras, {\it Invent.\ math.} {\bf 122} (1995), 341--357.

\item{18.} {\smcap T.\ Gannon, Ph.\ Ruelle} and {\smcap M.\ A.\ Walton},  
Automorphism modular invariants of current algebras, {\it Commun.\ Math.\ 
Phys.} {\bf 179} (1996), 121--156.

\item{19.} {\smcap G.\ Georgiev} and {\smcap O.\ Mathieu}, Cat\'egorie de
fusion pour les groupes de Chevalley, {\it C.\ R.\ Acad.\ Sci.\ Paris} {\bf 315}
 (1992), 659--662.

\item{20.} {\smcap F.\ M.\ Goodman} and {\smcap H.\ Wenzl},
Littlewood-Richardson coefficients for Hecke algebras at roots of unity,
{\it  Adv.\ Math.} {\bf 82} (1990), 244--265.

\item{21.} {\smcap Y.-Z.\ Huang} and {\smcap J.\ Lepowsky}, Intertwining
operator algebras and vertex tensor categories for affine Lie algebras,
{\it Duke Math.\ J.} {\bf 99} (1999), 113--134.

\item{22.} {\smcap V.\ G.\ Kac}, Simple Lie groups and the Legendre
symbol, in: ``Lie algebras, group theory, and partially ordered
algebraic structures'',  Lecture Notes in Math, Vol.\ 848,
Springer-Verlag, Berlin, 1981.
                                       
\item{23.} {\smcap V.\ G.\ Kac,} ``Infinite Dimensional Lie algebras'',
3rd edition, Cambridge University Press, Cambridge, 1990.

\item{24.} {\smcap V.\ G.\ Kac} and {\smcap D.\ H.\ Peterson},
Infinite-dimensional Lie algebras, theta
functions and modular forms, {\it Adv.\ Math.} {\bf 53} (1984), 125--264.

\item{25.} {\smcap V.\ G.\ Kac} and {\smcap M.\ Wakimoto}, Modular and conformal 
constraints in representation theory of affine algebras, {\it Adv.\ Math.} 
{\bf 70} (1988), 156--236.

\item{26.} {\smcap G.\ Lusztig}, Exotic Fourier transform, {\it Duke
Math.\ J.} {\bf 73} (1994), 227--241.

\item{27.} {\smcap I.\ G.\ Macdonald}, ``Symmetric functions and Hall
polynomials'', 2nd edition, Oxford University Press, New York, 1995.

\item{28.} {\smcap W.\ G.\ McKay, R.\ V.\ Moody} and {\smcap J.\ Patera},
Decomposition of tensor products of $E_8$ representations, {\it Alg., Groups
Geom.} {\bf 3} (1986), 286--328.

\item{29.} {\smcap W.\ G.\ McKay, J.\ Patera} and {\smcap D.\ W.\ Rand}, 
``Tables
of representations of simple Lie algebras'', Vol.\ 1, Centre de Recherches
Math\'ematiques, Univ\'ersit\'e de Montr\'eal, 1990.

\item{30.} {\smcap W.\ Nahm}, Conformal field theories, dilogarithms, and
3-dimensional manifolds, in: ``Interface between physics and mathematics'',
World-Scientific, 1994, (W.\ Nahm and J.-M.\ Shen, Eds.),
World-Scientific, Singapore,  1994.

\item{31.} {\smcap A.\ N.\ Schellekens}, Fusion rule automorphisms from
integer spin simple currents, {\it Phys.\ Lett.} {\bf B244} (1990), 255--260.

\item{32.} {\smcap V.\ G.\ Turaev}, ``Quantum invariants of knots and
3-manifolds'', Walter de Gruyter, Berlin, 1994.

\item{33.} {\smcap E.\ Verlinde}, Fusion rules and modular transformations
in 2D conformal field theory, {\it Nucl.\ Phys.} {\bf 300} (1988), 360--376.

\item{34.} {\smcap D.\ Verstegen}, New exceptional modular invariant
partition functions for simple Kac--Moody algebras, {\it Nucl.\ Phys.} {\bf B346}
(1990), 349--386.

\item{35.} {\smcap M.\ A.\ Walton}, Algorithm for WZW fusion rules: a 
proof, {\it Phys.\ Lett.} {\bf B241} (1990), 365--368.

\item{36.} {\smcap A.\ J.\ Wassermann}, Operator algebras and conformal
field theory, in: ``Proc.\ ICM, Zurich'', Birkh\"auser, Basel, 1995.

\item{37.} {\smcap E.\ Witten}, The Verlinde formula and the cohomology of
the Grassmannian, in: ``Geometry, Topology and Physics'', 
 International Press, Cambridge, MA, 1995.

\end